\documentclass[MSNbibl,number,citesort,dvips]{arxbj}
\usepackage{upgreek}
\usepackage{graphicx}

\aid{0}
\volume{20}
\issue{3}
\pubyear{2014}
\firstpage{1454}
\lastpage{1483}
\doi{10.3150/13-BEJ529} 

\makeatletter
\def\pi{\uppi}

\newcommand{\rrVert}{\Vert}
\newcommand{\rrvert}{\vert}
\newcommand{\llVert}{\Vert}
\newcommand{\llvert}{\vert}
\def\cal{\mathcal}
\newtheorem{theorem}{Theorem}
\newtheorem{lemma}{Lemma}
\newremark{rem}{Remark}
\newproclaim{definition}{Definition}
\newproclaim{assumption}{Assumption}
\newremark{example}{Example}
\makeatother

\begin{document}
\begin{frontmatter}

\title{Sojourn measures of Student and Fisher--Snedecor random fields}

\runtitle{Sojourn measures of Student and Fisher--Snedecor random fields}

\begin{aug}
\author[a]{\inits{N.}\fnms{Nikolai} \snm{Leonenko}\thanksref{a}\ead[label=e1]{LeonenkoN@cardiff.ac.uk}}
\and
\author[b]{\inits{A.}\fnms{Andriy} \snm{Olenko}\corref{}\thanksref{b}\ead[label=e2]{a.olenko@latrobe.edu.au}}
\runauthor{N. Leonenko and A. Olenko} 
\address[a]{School of Mathematics, Cardiff University, Senghennydd
Road, Cardiff CF24 4AG, United Kingdom. \printead{e1}}
\address[b]{Department of Mathematics and Statistics, La Trobe
University, Victoria, 3086, Australia.\\ \printead{e2}}

\end{aug}

\received{\smonth{6} \syear{2012}}
\revised{\smonth{12} \syear{2012}}

%
\begin{abstract}
Limit theorems for the volumes of excursion sets of weakly and strongly
dependent heavy-tailed random fields are proved. Some generalizations
to sojourn measures above moving levels and for cross-correlated
scenarios are presented. Special attention is paid to Student and
Fisher--Snedecor random fields. Some simulation results are also presented.
\end{abstract}

%

\end{frontmatter}

\section{Introduction}

Geometric characteristics of random surfaces play a crucial role in
areas such as geoscience, environmetrics, astrophysics, and medical
imaging, just to mention a few examples. Numerous real data have been
modelled as Gaussian random processes or fields and studying of their
excursion sets is now a well developed subject. Sojourn measures
provide a classical approach to addressing various applied problems
within this framework. There is a very rich literature on the topic,
therefore below we cite only some key publications related to our
approach. Good introductory references to some applications can be
found in \cite{adl,aza,bul,mar,nov}.

Sojourn measures of stochastic processes were studied extensively in a
number of contexts and explicit formulae for their statistical
characteristics were obtained for various scenarios, see, for example,
\cite{bra,kra,lin}, results for Gaussian stochastic processes with
long range dependence in~\cite{ber0,ber}, and also numerous references therein.
Unfortunately, one cannot expect that the same will occur for the
multidimensional situation. For random fields explicit formulae for the
excursion distributions are rarely known, see \cite{adl,bor}. Most
published papers concern only first two moments of sojourn measures.
However, it turned out that there are some interesting asymptotic
results in this area. Such results are usually the main tools for
statistical applications. It is natural to consider the volume of
excursion sets in a bounded observation window and to study its limit
behaviour as the window size grows. Some progress in this direction has
been made in \cite{adl0,bul,leo,leo0,liu1,liu2,mes}.

The approach taken in the paper continues this line of investigations.
The paper \cite{bul} studied central limit theorems for the volumes of
excursion sets of stationary quasi-associated random fields and
suggested two open problems: the extension of the results to different
classes of
random fields and the investigation of asymptotics for strongly
dependent structures.

\begin{figure}

\includegraphics{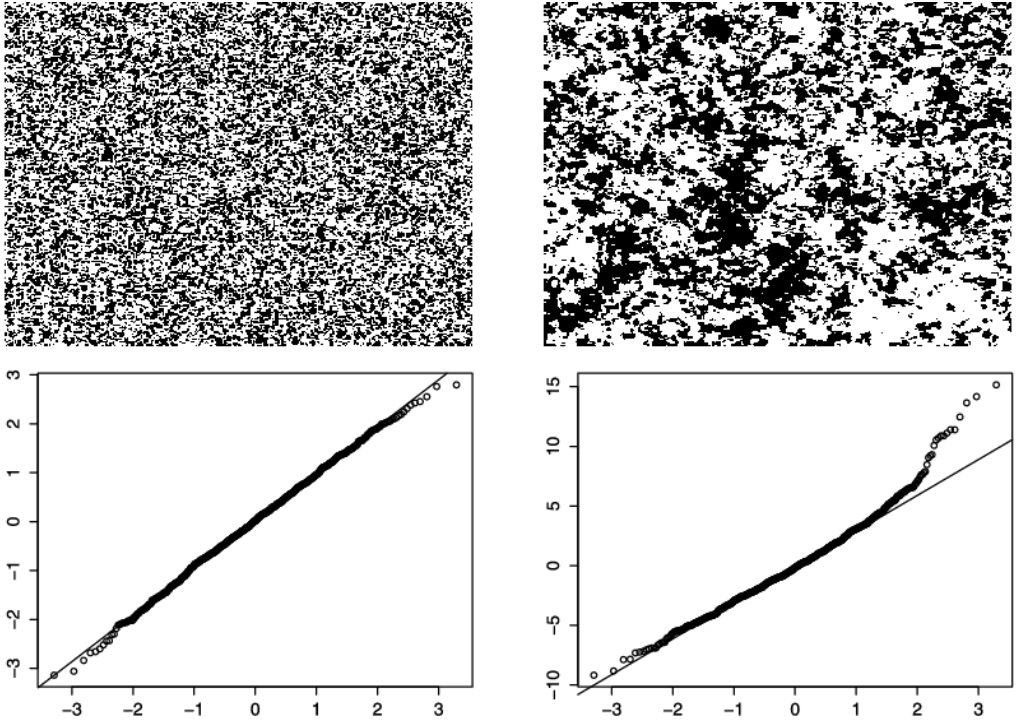}

\caption{Two-dimensional excursion sets and normal Q--Q plots of their
areas. The columns correspond
to short-range and long-range dependent models (from left to
right).}\label{fig1}
\end{figure}

In example Figure~\ref{fig1} the first row shows two-dimensional
excursion sets for realizations of two types of random fields (from
left to right): short-range dependent normal scale mixture model and
long-range dependent Cauchy model, consults Section~\ref{sec9}. The
excursion sets are shown in black colour. The Q--Q plots in the second
row, which correspond to the models shown above, suggest that the limit
law of the short-range dependent model is normal, while for the
long-range dependent model the data are not normally distributed.
Additional details about Figure~\ref{fig1} are provided in
Section~\ref{sec9}.\looseness=1\vadjust{\eject}


The paper has three aims. One is to provide explicit, albeit
asymptotic, formulae for the distribution of the volume of excursion
sets of a class of strongly dependent random fields. The second one is
to derive asymptotic results for heavy-tailed random fields. Finally,
the third aim is to generalize the previous findings to sojourn
measures above moving levels and for cross-correlated scenarios.

There is, therefore, a need for models that are able to display
strongly dependent heavy-tailed behaviour and yet are sufficiently
simple to allow
analysis. To obtain explicit results we detail the underlying structure
of random fields. Namely, a basic assumption of the analysis is that we
examine functionals of vector Gaussian random fields, in particular,
Student and Fisher--Snedecor random fields. Consult \cite
{ahm,cao1,cao2,wor} on excursion sets of chi-square, Student and
Fisher--Snedecor random fields and their importance for image analysis
and studies of brain function. Other results on sojourn measures of
chi-square random fields can be found in \cite{iv,leo,leo0,leo1}.

Minkowski functionals are widely used to characterise geometric
properties of random fields, in particular in the analysis of cosmic
microwave background radiation, see \cite{mar,nov}. In this paper we
investigate the first Minkowski functional of random fields and its
expansions into multidimensional Hermite polynomials, see some
one-dimensional/discrete counterparts in \cite{deh,dou}. To have a
complete account of results on asymptotic distributions of sojourn
measures for functions of vector random fields, we also prove
corresponding theorems for weakly dependent scenarios.

The remainder of the paper is structured as follows. In Sections~\ref{sec1}--\ref{sec3}, we introduce the necessary background from the
theory of random fields and briefly review some definitions and
notation on the first Minkowski functional, multidimensional Hermite
expansions, and Student and Fisher--Snedecor random fields. We start
Sections~\ref{sec4} and \ref{sec6} with generalizations and corrections
of some classical asymptotic results to arbitrary sets and vector
fields. With this in hand, we continue Sections~\ref{sec4} and~\ref
{sec6} by new results for the first Minkowski functional of Student and
Fisher--Snedecor random fields. In Section~\ref{sec6}, we also show how
to lift these results to sojourn measures above moving levels and for
cross-correlated underlying vector fields. Sections~\ref{proofs1} and~\ref{proofs2} provide the proofs of all theorems and lemmata in the
article. Simulation results on the limit distributions of areas of
excursion sets for two types of images are given in Section~\ref{sec9}.
Short conclusions are made in Section~\ref{sec10}.

In this paper, we only consider real-valued random fields. $\llvert
\cdot\rrvert  $ and $\llVert  \cdot\rrVert  $ denote the Lebesgue
measure and the distance in $\mathbb{R}^{d}$, respectively. In what
follows, we use the symbol $C$ to denote constants which are not
important for our discussion. Moreover, the same symbol $C$ may be used
for different constants appearing in the same proof.

\section{First Minkowski functional}\label{sec1}
In this section, we review the definition of the first Minkowski
functional and its relevant properties. More information about
stochastic Minkowski functionals and their links with the expected
Euler characteristics of excursion sets can be found in \cite{adl}.

We consider a measurable mean square continuous homogeneous isotropic
random field $S(x), x\in\mathbb{R}^{d}$, (see \cite{iv,leo1}) with
$\mathbf{E}S(x)=m$,
and the covariance function
\[
\mathrm{B}(r):=\operatorname{\mathbf{Cov}} \bigl( S(x),S(y) \bigr) =\int_{0}^{\infty
}Y_{d}(rz)
\,\mathrm{d}\Phi(z),\qquad x,y\in\mathbb{R}^{d},
\]
where $r:= \llVert  x-y\rrVert $, $\Phi(\cdot)$ is the isotropic
spectral measure, $Y_{d}(\cdot)$ is the spherical Bessel
function given by
\begin{eqnarray*}
Y_{1}(z)&:=&\cos z,
\\
Y_{n}(z)&:=&2^{(n-2)/2}\Gamma \biggl(\frac{n}{2} \biggr)
J_{(n-2)/2}(z) z^{(2-n)/2},\qquad z\geq0, n\geq2,
\end{eqnarray*}
$J_{\nu}(\cdot)$ is the Bessel function of the first kind of order $\nu
>-1/2$.

We define the marginal c.d.f. $H(\cdot)$ and p.d.f. $h(\cdot)$ of the field
$S(x)$ as follows:
\[
H(u)=\mathbf{P} \bigl\{ S(x)\leq u \bigr\} ,\qquad H(u)=\int_{-\infty}^u
h(z) \,\mathrm{d}z,\qquad u\in\mathbb{R}.
\]

\begin{definition} $S (x)$, $x\in\mathbb{R}^{d}$, is a homogeneous
isotropic random field possessing an absolutely continuous spectrum, if
there exists a function $f(\cdot)$ such that
\[
\Phi(z)=2\pi^{d/2}\Gamma^{-1} (d/2 )\int_0^z
u^{d-1}f(u)\,\mathrm{d}u, \qquad u^{d-1}f(u)\in\mathbf{L}_{1}(
\mathbb{R}_{+}).
\]
The function $f(\cdot)$ is
called the isotropic spectral density function of the field $S (x)$.
\end{definition}

Consider a Jordan-measurable convex bounded set $\Delta\subset\mathbb
{R}^{d}$, such that
$\llvert  \Delta\rrvert  >0$ and $\Delta$ contains the origin in its
interior. Let $\Delta(r),r>0$, be the homothetic image of
the set $\Delta$, with the centre of homothety in the origin and the
coefficient $r>0$,
that is, $\llvert  \Delta(r)\rrvert  =r^{d}\llvert  \Delta\rrvert $.

\begin{definition} The first Minkowski functional is defined as
\[
M_{r} \{ S \} :=\bigl\llvert \bigl\{x\in\Delta(r)\dvt S(x)>a(r)\bigr
\}\bigr\rrvert =\int_{\Delta(r)}\chi \bigl( S(x)>a(r) \bigr) \,\mathrm{d}x,
\]
where $\chi
(\cdot)$ is an indicator function and $a(r)$ is a continuous
non-decreasing function.
\end{definition}
In the
simplest case $a(r)=a$ is a constant. The functional $M_{r} \{
S \}$ has an interpretation
of the sojourn measure of the random field $S(x)$ above the constant
level $a$, or the moving level $a(r)$.

For the first Minkowski functional $M_{r} \{ S \}$ we obtain:
%
%
\begin{equation}
\label{exp} \mathbf{E} M_{r} \{ S \} =\llvert \Delta\rrvert
r^{d}\mathbf{P} \bigl\{ S(x)>a(r) \bigr\} =\llvert \Delta\rrvert
r^{d}\bigl(1-H\bigl(a(r)\bigr)\bigr)
\end{equation}
and
\[
\operatorname{\mathbf{Var}}M_{r} \{ S \} =\int_{\Delta(r)}\int
_{\Delta
(r)}\mathbf{P} \bigl\{ S(x)>a(r),S(y)>a(r) \bigr\} \,\mathrm{d}x\,\mathrm{d}y-
\bigl[ \mathbf {E}M_{r} \{ S \} \bigr] ^{2},
\]
or
\[
\operatorname{\mathbf{Var}}M_{r} \{ S \} =\int_{\Delta(r)}\int
_{\Delta(r)} \operatorname{\mathbf{Cov}} \bigl( \zeta(x),\zeta(y) \bigr) \,\mathrm{d}x\,\mathrm{d}y,
\]
where
$\zeta(x):=\chi (
S(x)>a(r) )$, $x\in\mathbb{R}^{d}$.
Therefore, it is important to investigate the integrals
\[
\int_{\Delta(r)}\int_{\Delta(r)}G\bigl(\llVert x-y
\rrVert \bigr)\,\mathrm{d}x \,\mathrm{d}y
\]
of various integrable Borel functions $G(\cdot)$.

Consider the uniform distribution on $\Delta(r)$
with the p.d.f. given by
\[
q_{\Delta(r)}(x)= \cases{\displaystyle\frac{1}{r^{d}\llvert  \Delta\rrvert  }, & \quad $\mbox{if } x\in\Delta (r)$;
\vspace*{2pt}
\cr
0, &\quad $\mbox{if } x\notin\Delta(r).$}
\]
Let $U$ and $V$ be two independent and uniformly distributed inside the
set $
\Delta(r)$ random vectors. We denote by $\psi_{\Delta(r)}(\rho
)$, $\rho\geq0$, the p.d.f. of the distance $\llVert  U-V\rrVert  $
between $U$ and $V$. Note that $\psi_{\Delta(r)}(\rho
)=0$ if $\rho> \mathrm{diam} \{ \Delta(r)
\}$.
Using the above notation, we obtain the representation
\begin{eqnarray}\label{dint}
\int_{\Delta(r)}\int_{\Delta(r)}G\bigl(\llVert x-y
\rrVert \bigr) \,\mathrm{d}x \,\mathrm{d}y&=& \llvert \Delta\rrvert ^{2}r^{2d}
\mathbf{E} G\bigl(\llVert U-V\rrVert \bigr)
\nonumber
\\[-8pt]
\\[-8pt]
\nonumber
 &=&\llvert \Delta\rrvert ^{2}r^{2d}\int
_{0}^{\mathrm{diam} \{ \Delta
(r) \} }G(\rho) \psi_{\Delta(r)}(\rho)\,\mathrm{d}\rho.
\end{eqnarray}

\begin{example} If $\Delta(r)$ is the ball $v(r):=\{x\in\mathbb{R
}^{d}\dvt \llVert  x\rrVert  <r\}$ then
\[
\psi_{v(r)}(\rho)=d \rho^{d-1}r^{-d}I_{1-(\rho/2r)^{2}}
\biggl( \frac{d+1}{2} ,\frac{1}{2} \biggr) ,\qquad 0\leq\rho\leq2r,
\]
where
%
%
\begin{equation}
\label{beta} I_{\mu}(p,q):=\frac{\Gamma(p+q)}{\Gamma(p) \Gamma(q)}\int_{0}^{
\mu}t^{p-1}(1-t)^{q-1}\,\mathrm{d}t,\qquad
\mu\in( 0,1], p>0, q>0 ,
\end{equation}
is the incomplete beta function, see \cite{iv}.

Several expressions for $\psi_{v(r)}(\rho)$, $0\leq\rho\leq2r$,
are given below:
\begin{enumerate}[$d=1$:]
\item[$d=1$:] $\psi_{v(r)}(\rho)=\frac{1}{r} ( 1-\frac{\rho
}{2r} )$,
\item[$d=2$:] $\psi_{v(r)}(\rho)=\frac{4\rho}{\pi r^2} (
\arccos\frac{\rho}{2r}-\frac{\rho}{2r}\sqrt{1- ( \frac{\rho}{2r}
) ^{2}} )$,
\item[$d=3$:] $\psi_{v(r)}(\rho)= \frac{3\rho^{2} }{r^{3}}  (
1-\frac{\rho}{2r} )
^{2} ( 1+\frac{\rho}{4r} )$.\vadjust{\goodbreak}
\end{enumerate}

If one considers the functional
\[
F_r(\zeta)=\int_{v(r)}\zeta(x) \,\mathrm{d}x,
\]
then
\begin{eqnarray*}
\operatorname{\mathbf{Var}} F_r(\zeta)&=&\int_{v(r)}\int
_{v(r)}\tilde{\mathrm{B}}\bigl(\llVert x-y\rrVert \bigr) \,\mathrm{d}x \,\mathrm{d}y
=\bigl\llvert v(1)\bigr\rrvert ^{2}r^{2d}\mathbf{E} \tilde{
\mathrm{B}}\bigl(\llVert U-V\rrVert \bigr)
\\
&=&\frac{4\pi^{d}}{d \Gamma^{2}(d/2)}r^{d}\int_{0}^{2r}z^{d-1}
\tilde {\mathrm{B}}(z) I_{1-({z}/{2r})^{2}} \biggl( \frac{d+1}{2},
\frac{1}{2} \biggr) \,\mathrm{d}z,
\end{eqnarray*}
where $\tilde{\mathrm{B}}(\cdot)$ is a covariance function of $\zeta(x)$.
\end{example}

For some random fields these formulae can be specified,
however the asymptotic analy\-sis is difficult. Therefore, we will use
an approach based on multidimensional Hermite expansions.

\section{Multidimensional Hermite expansions}\label{sec2}

Let $H_{k}(u)$, $k\geq0$, $u\in\mathbb{R}$, be the Hermite polynomials, see \cite{pec}.
%
%

\begin{lemma}\cite{pec}\label{lem1} Let $(\xi_{1},\ldots,\xi
_{2p})$ be $2p$-dimensional
zero mean Gaussian vector with
\begin{eqnarray*}
\mathbf{E}\xi_{j}\xi_{k} = \cases{ 1, &\quad$\mbox{if } k=j$;
\vspace*{2pt}
\cr
r_{j}, &\quad $\mbox{if } k=j+p \mbox{ and } 1\leq j\leq
p;$\vspace*{2pt}
\cr
0, &\quad $\mbox{otherwise.}$}
\end{eqnarray*}
Then
\[
\mathbf{E} \prod_{j=1}^{p}H_{k_{j}}(
\xi_{j})H_{m_{j}}(\xi _{j+p})=\prod
_{j=1}^{p}\delta_{k_{j}}^{m_{j}}
k_{j}! r_{j}^{k_{j}}.
\]
\end{lemma}

Let us denote
\[
e_{\nu}(w):=\prod_{j=1}^{p}H_{k_{j}}(w_{j}),
\]
where $w=(w_{1},\ldots
,w_{p})^{\prime}\in\mathbb{R}^{p}$, $\nu=(k_{1},\ldots,k_{p})\in
\mathbb{Z}^{p}$, and all $k_{j}\geq0$ for $j=1,\ldots,p$.

The summation theorem for Hermite polynomials \cite{gra}, formula (8.958.1)
states that
%
%
\begin{equation}
\label{sumH} H_k \biggl(\frac{\sum_{j=1}^pa_jw_j}{\sum_{j=1}^pa^2_j} \biggr)=
\frac{k!}{ (\sum_{j=1}^pa^2_j
)^{k/2}}\sum_{k_1+\cdots+k_p=k} \prod
_{j=1}^p\frac{a_j^{k_j}}{k_j!}H_{k_j}(w_j).
\end{equation}

The polynomials $\{e_{\nu}(w)\}_{\nu}$ form a complete orthogonal system
in the Hilbert space
\begin{eqnarray*}
\mathbf{L}_{2}\bigl(\mathbb{R}^{p},\phi\bigl(\llVert w
\rrVert \bigr) \,\mathrm{d}w\bigr) &=& \biggl\{ G\dvt \int_{
\mathbb{R}^{p}}G^{2}(w)
\phi\bigl(\llVert w\rrVert \bigr) \,\mathrm{d}w<\infty \biggr\},
\\
\phi\bigl(\llVert w\rrVert \bigr) &=&\prod_{j=1}^{p}
\phi(w_{j}),\qquad \phi (w_{j})=\frac{1}{\sqrt{2\pi}}\mathrm{e}^{-{w_{j}^{2}}/{2}}.
\end{eqnarray*}

An arbitrary function $G(w)\in\mathbf{L}_{2}(\mathbb{R}^{p},\phi(\llVert
w\rrVert  )  \,\mathrm{d}w)$ admits the mean-square convergent expansion
%
%
\begin{equation}
\label{herm} G(w)=\sum_{k=0}^{\infty}\sum
_{\nu\in N_{k}}\frac{C_{\nu}e_{\nu}(w)
}{\nu!},\qquad  C_{\nu}:=\int
_{\mathbb{R}^{p}}G(w)e_{\nu
}(w)\phi\bigl(\llVert w\rrVert
\bigr) \,\mathrm{d}w,
\end{equation}
where $
N_{k} :=\{(k_{1},\ldots,k_{p})\in\mathbb{Z}^{p}\dvt \sum_{j=1}^{p}k_{j}=k, \mbox{ all } k_{j}\geq0 \mbox{ for } j=1,\ldots,p\},
$ $\nu! :=k_{1}!\cdots k_{p}!$

By Parseval's identity
%
%
\begin{equation}
\label{par} \sum_{k=0}^\infty\sum
_{\nu\in N_{k}}\frac{C_{\nu}^{2}}{\nu!} =\int_{R^{p}}G^2(w)
\phi\bigl(\llVert w\rrVert \bigr) \,\mathrm{d}w.
\end{equation}

\begin{definition} Let $G(w)\in\mathbf{L}_{2}(\mathbb{R}^{p},\phi
(\llVert  w\rrVert  )
\,\mathrm{d}w)$ and there exist an integer $\kappa\geq1$ such that $C_{\nu}=0$,
for all $\nu
\in N_{k}$, $0\leq k\leq\kappa-1$, but $C_{\nu}\neq0$ for at least one
tuple $\nu=(k_{1},\ldots,k_{p})\in N_{\kappa}$. Then $\kappa$ is
called the Hermite rank of $G(\cdot)$ and denoted by
$H \operatorname{rank} G$.
\end{definition}

Let $\eta(x)=[\eta_{1}(x),\ldots,\eta_{p}(x)]^{\prime}$, $x\in
\mathbb{R}
^{d}$, be a measurable mean-square continuous homogeneous isotropic
vector Gaussian random field, see Section~5 in \cite{leo1}, Section~1.2.
Suppose that the components $\eta_{1}(\cdot),\ldots,\eta_{p}(\cdot)$
are independent, $\mathbf{E}\eta_j(0)=0$, $\mathbf{E}\eta_j^2 (0)=1$,
and $\mathbf{E}\eta_j (0)\eta_j (x)=\mathrm{B}_{jj}(\llVert  x\rrVert
)$, $1\leq j\leq p$.

If $G(w)\in\mathbf{L}_{2}(\mathbb{R}^{p},\phi
(\llVert  w\rrVert  )  \,\mathrm{d}w)$ then the integral functional $F(\eta)=\int_{\Delta(r)}G(\eta(x))  \,\mathrm{d}x$
can be represented as
\[
F(\eta)=\sum_{k=0}^{\infty}\sum
_{\nu\in N_{k}}\frac{C_{\nu}}{\nu!} \int_{\Delta(r)}e_{\nu}
\bigl(\eta(x)\bigr) \,\mathrm{d}x.
\]
Therefore the expectation of $F(\eta)$ is
%
%
\begin{equation}
\label{expg}\mathbf{E}F(\eta)=\bigl|\Delta(r)\bigr|C_{(0,\ldots,0)},
\end{equation}
while by Lemma~\ref{lem1} the variance is equal
%
%
\begin{equation}
\label{varg} \operatorname{\mathbf{Var}} F(\eta)= \sum_{k=0}^{\infty}
\sum_{\nu\in N_{k}}\frac{C_{\nu}^{2}}{\nu!} \int
_{\Delta(r)}\int_{\Delta(r)}\prod
_{j=1}^{p}\mathrm {B}_{jj}^{k_{j}}
\bigl(\llVert x-y\rrVert \bigr)\,\mathrm{d}x\,\mathrm{d}y.
\end{equation}

\section{Student and Fisher--Snedecor random fields}\label{sec3}
In this section, we introduce two main models investigated in the
paper, namely, Student and Fisher--Snedecor random fields proposed for
studies of brain
function in~\cite{wor}.

Let us consider the vector random field
\[
\eta(x)=\bigl[ \eta_{1}(x),\ldots,\eta_{m}(x),
\eta_{m+1}(x),\ldots,\eta_{m+n}(x)\bigr]',
\]
which consists of $n+m$ independent copies of a measurable mean-square
continuous homogeneous isotropic zero-mean and unit variance Gaussian
random field $
\eta_{1}(x)$, $x\in\mathbb{R}^{d}$.

\begin{definition}\label{def3} The Fisher--Snedecor random field
$F_{m,n}(x)$, $x\in\mathbb{R}^{d}$, is defined by
\[
F_{m,n}(x):=\frac{({1}/{m}) ( \eta_{1}^{2}(x)+\cdots+\eta
_{m}^{2}(x) )
}{({1}/{n})(\eta_{m+1}^{2}(x)+\cdots+\eta_{m+n}^{2}(x))},\qquad  x\in \mathbb{R}^{d}.
\]
\end{definition}

The random field $ F_{m,n}(x)$, $x\in\mathbb{R}^{d}, $ has the
marginal Fisher--Snedecor distribution with the p.d.f.
\[
h(u)=\frac{m^{{m}/{2}}n^{{n}/{2}}\Gamma ({(m+n)}/{2}
)}{\Gamma ({m}/{2} )\Gamma ({n}/{2} )}\cdot \frac{u^{{m}/{2}-1}}{(n+mu)^{{(n+m)}/{2}}}, \qquad u\in \lbrack0,\infty),
\]
and the c.d.f.
%
%
\begin{equation}
\label{Fmn} H(u)=I_{{mu}/{(n+mu)}} \biggl(\frac{m}{2},\frac{n}{2}
\biggr).
\end{equation}

By properties of the Fisher--Snedecor distribution
%
%
\[
\mathbf{E} \bigl[ F_{m,n}(x) \bigr] ^{r}=\frac{\Gamma (
{(m+2k)}/{2} )\Gamma ({(n-2k)}/{2} )}{\Gamma (
{m}/{2} )\Gamma ({n}/{2} )}
\biggl(\frac{n}{m} \biggr)^{r},\qquad n>2r.
\]

\begin{definition}\label{def4} The Student random field $
T_{n}(x)$, $x\in\mathbb{R}^{d}$, is defined by
\[
T_{n}(x):=\frac{\eta_{1}(x)}{\sqrt{({1}/{n}) ( \eta
_{2}^{2}(x)+\cdots
+\eta_{n+1}^{2}(x) ) }},\qquad x\in\mathbb{R}^{d}.
\]
\end{definition}
It has the marginal Student $t_{n}$-distribution with the p.d.f.
\[
h(u)=\frac{\Gamma ( {(n+1)}/{2} ) }{\sqrt{n\pi}\Gamma (
{n}/{2} )}\cdot \biggl( 1+\frac{u^{2}}{n} \biggr)^{-{(n+1)}/{2}}
,\qquad u\in\mathbb{R},
\]
and the c.d.f.
%
%
\begin{equation}
\label{SH}H(u)=\frac12+\frac12 \biggl(1-I_{
{n}/{(n+u^2)}} \biggl(
\frac{n}{2},\frac{1}{2} \biggr) \biggr)\cdot\operatorname{sgn}(u),
\end{equation}
where $\operatorname{sgn}(\cdot)$ is the signum function.

The $r$th moments of $T_{n}(x)$ exist when $n>r$ and for $k\in\mathbb
N$ we have
\[
\mathbf{E} \bigl\{ T_{n}(x) \bigr\} ^{r}= \cases{0, &\quad$
\mbox{if } r=2k-1<n$; \vspace*{2pt}
\cr
\displaystyle\frac{\Gamma ({(r+1)}/{2} )
\Gamma ({(n-r)}/{2} )n^{r/2}}{\sqrt{\pi}\Gamma (
{n}/{2} )}, &\quad $\mbox{if } r=2k<n$.
}
\]

Note that $ [ T_{n}(x) ] ^{2}=F_{1,n}(x), x\in\mathbb{R}^{d}$.

\begin{rem}The right-hand tail of the p.d.f. of the $F_{m,n}$-distribution
decreases as $x^{-{(n+2)}/{2}}$.
The left and the right-hand tails of the p.d.f. of the $t$-distribution
decrease as $\llvert  x\rrvert ^{-n-1}$. Thus, both Student and
Fisher--Snedecor random fields have
heavy-tailed marginal distributions.
\end{rem}

\section{Central limit theorem for functionals of weakly dependent
vector random fields}\label{sec4}

In this section we present some analogues of results in \cite{arc1,arc2,bre,har} for the case of integrals of weakly dependent vector
random fields. Then, we apply these results to Fisher--Snedecor and
Student random fields.

Let $\eta(x)=[\eta_{1}(x),\ldots,\eta_{p}(x)]^{\prime}$, $x\in
\mathbb{R}
^{d}$, be a measurable mean-square continuous homogeneous isotropic
vector Gaussian random field with $\mathbf{E}\eta
(x)=0$ and covariance matrix
\[
\mathbf{B}\bigl(\llVert x\rrVert \bigr)=\mathbf{E}\eta(0)\eta(x)^{\prime
}=
\bigl(\mathrm{B}_{ij}\bigl(\llVert x\rrVert \bigr)\bigr)_{1\leq
i,j\leq p}.
\]

First, we need an auxiliary statement which is similar to Theorem~1 in
\cite{bre}.
Let $\Box_{u,v}(r):=\{x\in\mathbb{R}^d\dvt ru_i<x_i\le rv_i, i=1,\ldots,d\}
$, where $u,v\in\mathbb{R}^d$ and $u_i< v_i$ for all $1\le i\le d$.
We will use the notation
\[
\psi(x):=\max_{1\le i\le p}\sum_{j=1}^p
\bigl|\mathrm{B}_{ij}\bigl(\llVert x\rrVert \bigr)\bigr|.
\]

\begin{lemma}\label{th1} Suppose that the function $G(\cdot)$ has
Hermite rank $\kappa\ge1$, the covariance matrix of the vector field
$\eta(x)$ satisfies the conditions $\psi(\cdot)\in\mathbf{L}_\kappa
(\mathbb{R}^d)$ and $\psi(x)\le1$ for all $x\in\mathbb{R}^d$, and
\[
\sigma^2:=\int_{\mathbb R^d} \mathbf{E}\bigl[G\bigl(
\eta(0)\bigr)G\bigl(\eta(x)\bigr)\bigr]\,\mathrm{d}x\neq0.
\]
Then
\[
r^{-d/2}\int_{\Box_{u,v}(r)} G\bigl(\eta(x)\bigr)\,\mathrm{d}x\stackrel{
\mathcal{D}} {\to} Y_{\Box_{u,v}},\qquad r\to\infty,
\]
where $|\Box_{u,v}(1)|^{-1/2} Y_{\Box_{u,v}}\sim N(0,\sigma^2)$,
$|\Box_{u,v}(1)|=\prod_{i=1}^d(v_i-u_i)$.

If $\Box_{u^{(1)},v^{(1)}}(1) \cap\Box_{u^{(2)},v^{(2)}}(1)=\varnothing
$, $u^{(i)},v^{(i)}\in\mathbb{R}^d$, $i=1,2$, then the random variables
$Y_{\Box_{u^{(1)},v^{(1)}}}$ and $Y_{\Box_{u^{(2)},v^{(2)}}}$ are independent.
\end{lemma}
The proof of the lemma is based on Lemma~\ref{lem1}, the diagram
formula and ideas in \cite{bre}, see also~\cite{arc1,arc2} for vector
processes, and the application of the diagram technique for random
fields in~\cite{iv}. The assumption $\psi(\cdot)\in\mathbf{L}_\kappa
(\mathbb{R}^d)$ can be weakened, consult, for example, the conditions~(1.4$'$) and (1.4$''$) in Theorem~1$'$
\cite{bre}. The most recent results can be found in \cite{bar,iv1,nua,pec}.

The following result generalizes Theorem~4 in \cite{arc1} to the case
of integrals of weakly dependent vector random fields.

\begin{theorem}\label{th2} If the conditions of Lemma~\ref{th1}
are satisfied, then
\[
r^{-d/2}\int_{\Delta(r)} G\bigl(\eta(x)\bigr)\,\mathrm{d}x\stackrel{
\mathcal{D}} {\to} Y_{\Delta},\qquad r\to\infty,
\]
where $|\Delta|^{-1/2} Y_{\Delta}\sim N(0,\sigma^2)$.\vadjust{\goodbreak}
\end{theorem}

\begin{rem}
The central limit theorems for the volumes of excursion sets of
stationary quasi-associated random fields were proved in \cite
{bul,mes}. The approach used in the papers did not require the isotropy
of Gaussian fields. However, it was assumed that the continuous
covariance function is ${\cal O}(\llVert  x\rrVert ^{-\alpha})$, $\alpha
>d$, when $\llVert  x\rrVert \to\infty$. We obtain the central limit
theorems for homogeneous isotropic random fields but under different
conditions. Namely, it follows from (\ref{eq1}) that only the
integrability of the covariance functions is required.
\end{rem}


In the next two theorems we consider sojourn measures of
Fisher--Snedecor and Student random fields above the constant level
$a(r)\equiv a$. In the notation of Sections~\ref{sec1} and \ref{sec3},
for the Fisher--Snedecor random field $p=m+n$ and the first Minkowski
functional takes the form
\[
M_{r} \{ F_{m,n} \} =\bigl\llvert \bigl\{x\in\Delta(r)\dvt
F_{m,n}(x)>a\bigr\} \bigr\rrvert =\int_{\Delta(r)}\chi
\bigl( F_{m,n}(x)>a \bigr) \,\mathrm{d}x.
\]

\begin{theorem}\label{th3} If the covariance matrix of the
Fisher--Snedecor random field $ F_{m,n}(x)$, $x\in\mathbb{R}^d$,
satisfies the two conditions: $\sup_{x\in R^d}\psi(x)\le1$ and
$\psi(\cdot)\in\mathbf{L}_2({\mathbb R^d})$, then
\[
r^{-d/2}M_{r} \{ F_{m,n} \}-\llvert \Delta\rrvert
r^{d/2} \biggl(1-I_{{ma}/{(n+ma)}} \biggl(\frac{m}{2},
\frac{n}{2} \biggr) \biggr)\stackrel{\mathcal{D}} {\to} Y_{\Delta},\qquad r
\to\infty,
\]
where $|\Delta|^{-1/2} Y_{\Delta}\sim N(0,\sigma_{F}^2(a))$,
$I_{\mu}(p,q)$ is defined by (\ref{beta}),
\[
\sigma_{F}^2(a):=\int_{\mathbb R^d}
\mathbf{E} \bigl[\chi \bigl( F_{m,n}(0)>a \bigr)\chi \bigl(
F_{m,n}(x)>a \bigr) \bigr]\,\mathrm{d}x.
\]
\end{theorem}

For the Student, random field $p=n+1$ and the first Minkowski functional
for the constant level $a$ is
\[
M_{r} \{ T_{n} \} =\bigl\llvert \bigl\{x\in\Delta(r)\dvt
T_{n}(x)>a\bigr\}\bigr\rrvert =\int_{\Delta(r)}\chi
\bigl( T_{n}(x)>a \bigr) \,\mathrm{d}x.
\]
%
\begin{theorem}\label{thS} If the covariance matrix of the Student
random field $T_{n}(x)$, $x\in\mathbb{R}^d$, satisfies the two
conditions: $\sup_{x\in R^d}\psi(x)\le1$ and
$\psi(\cdot)\in\mathbf{L}_1({\mathbb R^d})$, then
\[
r^{-d/2}M_{r} \{ T_{n} \}-\llvert \Delta\rrvert
r^{d/2} \biggl(\frac12-\frac12 \biggl(1-I_{{n}/{(n+a^2)}} \biggl(
\frac{n}{2},\frac
{1}{2} \biggr) \biggr)\cdot\operatorname{sgn}(a)
\biggr)\stackrel{\mathcal{D}} {\to } \tilde{Y}_{\Delta},\qquad r\to\infty,
\]
where $|\Delta|^{-1/2} \tilde{Y}_{\Delta}\sim N(0,\sigma_T^2)$,
\[
\sigma_T^2:=\int_{\mathbb R^d} \mathbf{E}
\bigl[\chi \bigl( T_{n}(0)>a \bigr)\chi \bigl( T_{n}(x)>a
\bigr) \bigr]\,\mathrm{d}x.
\]
\end{theorem}

\section{Proofs of the results of Section~\texorpdfstring{\protect\ref{sec4}}{5}}\label{proofs1}

\begin{pf*}{Proof of Lemma~\ref{th1}} The lemma can be proved by a
modification of the proof of Theorem~1 \cite{bre}\vadjust{\goodbreak} using vector results
in \cite{arc1,arc2}. To avoid lengthy repetitions, we only state
required changes to Theorem~1 \cite{bre}.

The first step is the replacement of the function of a single variable
$H(t)$ in Theorem~1 by the function of multiple variables $G(x)$ and
use vector notation and conditions on the covariance matrix presented
in \cite{arc2}. Then, it is straightforward to replace the summation
over the sets $B(n,N):=\{s=(s_1,\ldots,s_d)\in\mathbb{Z}^d\dvt Nn_i<s_i\le N
(n_i+1), i=1,\ldots,d\}$, by the integration over the multidimensional
parallelepipeds $\Box_{u,v}(r):=\{x\in\mathbb{R}^d\dvt ru_i<x_i\le rv_i, i=1,\ldots,d\}$. Finally, using integrals instead of sums in Theorem~4
\cite{arc1} we obtain $\lim_{r\to\infty}r^{-d}\operatorname{\mathbf{Var}} (\int_{\Box_{u,v}(r)} G(\eta(x))\,\mathrm{d}x )$ and the expression for $\sigma^2$.

The condition $\psi(\cdot)\in\mathbf{L}_\kappa(\mathbb{R}^d)$
guarantees that cross-correlation functions of all components of $\eta
(x)$ are also in $\mathbf{L}_\kappa(\mathbb{R}^d)$.
\end{pf*}

\begin{pf*}{Proof of Theorem~\ref{th2}} Let us consider a coverage
of $\Delta(r)$ by the finite union $\Box_J(r):=\bigcup_{j\in J}\Box
_{u^{(j)},v^{(j)}}(r)$ of the disjoint multidimensional parallelepipeds
$\{\Box_{u^{(j)},v^{(j)}}(r),  j\in J\}$, with the following properties:
\begin{enumerate}
\item$\Box_J(r)$ is a decreasing nested sequence of sets when $r$ is
fixed and $|J|\to\infty$;
\item$\Delta\subset\Box_J(1)$;
\item$|\Box_J(1)\setminus\Delta|\to0$, when $|J|\to\infty$.
\end{enumerate}
The existence of such $\Box_J(1)$ follows form the fact that $\Delta$
is a Jordan-measurable set.

By Lemma~\ref{th1}, we obtain
\[
r^{-d/2}\int_{\Box_J(r)} G\bigl(\eta(x)\bigr)\,\mathrm{d}x\stackrel{
\mathcal{D}} {\to} Y_{\Box
_J}, \qquad r\to\infty,
\]
where $|\Box_J(1)|^{-1/2} Y_{\Box_J}\sim N(0,\sigma^2)$.

By the properties of $\Box_J(r)$, we get $Y_{\Box_J}\stackrel{\mathcal
{D}}{\to} Y_{\Delta}, |J|\to\infty$.

As $\psi(x)\le1$, then by Lemma~1 \cite{arc1}
%
%
\begin{equation}
\label{eq1} \bigl|\mathbf{E}\bigl[G\bigl(\eta(x)\bigr)G\bigl(\eta\bigl(x^{(1)}
\bigr)\bigr)\bigr]\bigr|\le\psi^\kappa \bigl(\bigl\llVert x-x^{(1)}
\bigr\rrVert \bigr)\mathbf{E}G^2\bigl(\eta(0)\bigr),\qquad
x,x^{(1)}\in\mathbb{R}^d.
\end{equation}

It follows from inequality (\ref{eq1}) that
\begin{eqnarray}\label{eq2}
&& r^{-d}\operatorname{\mathbf{Var}} \biggl( \int_{\Box_J(r)} G\bigl(
\eta(x)\bigr)\,\mathrm{d}x -\int_{\Delta
(r)} G\bigl(\eta(x)\bigr) \,\mathrm{d}x \biggr)
\nonumber\\
&&\quad= r^{-d} \int_{\Box_J(r)\setminus\Delta(r)}\int_{\Box_J(r)\setminus
\Delta(r)}
\mathbf{E} G\bigl(\eta(x)\bigr)G\bigl(\eta\bigl(x^{(1)}\bigr)\bigr) \,\mathrm{d}x
\,\mathrm{d}x^{(1)}
\nonumber
\\[-8pt]
\\[-8pt]
\nonumber
&&\quad\le\frac{\mathbf{E}G^2(\eta(0))}{r^{d}} \int_{\Box_J(r)\setminus\Delta
(r)}\int_{\Box_J(r)\setminus\Delta(r)}
\psi^\kappa \bigl(\bigl\llVert x-x^{(1)}\bigr\rrVert \bigr) \,\mathrm{d}x
\,\mathrm{d}x^{(1)}
\\
&&\quad\le\bigl|\Box_J(1)\setminus\Delta\bigr|\cdot\mathbf
{E}G^2\bigl(\eta(0)\bigr) \int_{\mathbb{R}^d}
\psi^\kappa \bigl(\llVert x\rrVert \bigr) \,\mathrm{d}x.\nonumber
\end{eqnarray}
Finally, by property~3 of $\Box_J(r)$ the upper bound in (\ref{eq2})
approaches 0 when $|J|\to\infty$, which completes the proof.
\end{pf*}

\begin{pf*}{Proof of Theorem~\ref{th3}}
Note that by (\ref{Fmn})
\[
\mathbf{E} \bigl(\chi \bigl( F_{m,n}(x)>a \bigr) \bigr)= \mathbf{P}
\bigl( F_{m,n}(x)>a \bigr)=1-I_{{ma}/{(n+ma)}} \biggl(
\frac{m}{2},\frac
{n}{2} \biggr).
\]
Then it follows from (\ref{exp}) that
\[
\mathbf{E} \biggl(r^{-d/2}\int_{\Delta(r)} \chi \bigl(
F_{m,n}(x)>a \bigr) \,\mathrm{d}x \biggr)=\llvert \Delta\rrvert r^{d/2}
\biggl(1-I_{{ma}/{(n+ma)}} \biggl(\frac{m}{2},\frac{n}{2} \biggr)
\biggr)
\]
and we obtain the following representation
\[
r^{-d/2}\int_{\Delta(r)} \bigl(\chi \bigl(
F_{m,n}(x)>a \bigr)-\mathbf {E} \bigl(\chi \bigl( F_{m,n}(x)>a
\bigr) \bigr) \bigr)\,\mathrm{d}x= r^{-d/2}\int_{\Delta(r)}G \bigl(
\eta(x) \bigr) \,\mathrm{d}x,
\]
where
%
%
\begin{equation}
\label{fG} G ( w )=\chi \biggl( \frac{({1}/{m}) (
w_{1}^{2}+\cdots+w_{m}^{2} ) }{({1}/{n}) (w_{m+1}^{2}+\cdots
+w_{m+n}^{2} )}>a \biggr)+I_{{ma}/{(n+ma)}}
\biggl(\frac{m}{2},\frac
{n}{2} \biggr)-1.
\end{equation}
$G ( \cdot )$ is a symmetric function with respect to the
origin. Hence, $C_{\nu}=0$ for all $\nu
\in N_{1}$. However, $C_{\nu}\neq0$ for such tuples $\nu
=(k_{1},\ldots,k_{m+n})\in N_{2}$ that exactly one $k_{i}=2$
(expressions for coefficients $C_{\nu}$, $\nu\in N_{2}$, will be
given in Theorem~\ref{th7}).

Therefore, $H \operatorname{rank} G=2$ and we can apply Theorem~\ref{th2}
which completes the proof.
\end{pf*}

\begin{pf*}{Proof of Theorem~\ref{thS}}
It is easy to obtain the statement of the theorem following steps
analogous to the proof of Theorem~\ref{th3}.

Using (\ref{SH}), we conclude that
\[
\mathbf{E} \bigl(\chi \bigl( T_{n}(x)>a \bigr) \bigr)= \frac12-\frac
12 \biggl(1-I_{{n}/{(n+a^2)}} \biggl(\frac{n}{2},\frac{1}{2} \biggr)
\biggr)\cdot\operatorname{sgn}(a).
\]

Therefore,
\[
r^{-d/2}\int_{\Delta(r)} \bigl(\chi \bigl(
T_{n}(x)>a \bigr)-\mathbf {E} \bigl(\chi \bigl(T_{n}(x)>a
\bigr) \bigr) \bigr)\,\mathrm{d}x= r^{-d/2}\int_{\Delta(r)}\tilde{G}
\bigl(\xi(x) \bigr) \,\mathrm{d}x,
\]
where
%
%
\begin{equation}
\label{tG} \tilde{G} ( w )=\chi \biggl( \frac{w_{1}}{\sqrt{({1}/{n})
(w_{2}^{2}+\cdots
+w_{n+1}^{2} )}}>a \biggr)+\frac12
\biggl(1-I_{{n}/{(n+a^2)}} \biggl(\frac{n}{2},\frac{1}{2} \biggr)
\biggr)\cdot\operatorname{sgn}(a)-\frac12.
\end{equation}
For $\tilde{G} ( \cdot )$ the coefficient $C_{(1,0,\ldots,0)}\neq
0$, $(1,0,\ldots,0) \in N_{1}$, (expressions for coefficients $C_{\nu}$,
$\nu\in N_{1}$, will be given in Theorem~\ref{th6}). Therefore, $H
\operatorname{rank} \tilde{G}=1$ and the application of Theorem~\ref{th2}
completes the proof.
\end{pf*}

\section{Non-central limit theorem for functionals of strongly
dependent vector random fields}\label{sec6}

In this section, we first present corrections and generalizations to
arbitrary sets of some results for random fields in \cite{iv}, Section~2.10,
\cite{leo1}, Sections~2.4 and 3.4, and \cite{leoole0}. Consult also the
pioneering papers \cite{dob,ta0,ta2} and the book \cite{ber} on
non-central limit theorems and the Hermite polynomials approach. In the
rest of this section, we apply the developed technique to
Fisher--Snedecor and Student random fields.

\begin{assumption}\label{ass1} { Let $\eta(x)=[\eta_{1}(x),\ldots
,\eta_{p}(x)]^{\prime}$, $x\in\mathbb{R}^{d}$, be a vector
homogeneous isotropic Gaussian random
field with $\mathbf{E}\eta(x)=0$ and covariance matrix
\[
\mathbf{\tilde{B}}(0)=\mathcal{I}, \qquad\mathbf{\tilde{B}}\bigl(\llVert x\rrVert
\bigr)=\mathbf{E}\eta(0) \eta(x)^{\prime}=\mathcal{I}\cdot\llVert x\rrVert
^{-\alpha}L\bigl(\llVert x\rrVert \bigr),\qquad \alpha>0,
\]
where $\mathcal{I}$ is the unit matrix of size $p$, $L(\llVert \cdot
\rrVert )$ is a function slowly varying at infinity.}
\end{assumption}

We investigate the random variables
\[
K_r :=\int_{\Delta(r)}G_r\bigl(\eta(x)
\bigr)\,\mathrm{d}x \quad\mbox{and}\quad K_{r,\kappa} :=\sum_{\nu\in N_\kappa}
\frac{C_\nu(r)}{\nu!} \int_{\Delta(r)}e_\nu\bigl(\eta(x)
\bigr)\,\mathrm{d}x,
\]
where $C_\nu(r)$ are coefficients of the Hermite series (\ref{herm})
of the function $G_r(\cdot)$ for fixed $r$.

\begin{theorem}\label{th4}
Suppose that $\eta(x)$ satisfies Assumption~\ref{ass1} for
$\alpha\in
(0,d/\kappa)$, for each sufficiently large $r$ $H \operatorname{rank}
G_{r}(\cdot)=\kappa\ge1$, and
%
%
\begin{equation}
\label{con1} \biggl( \sum_{\nu\in N_{\kappa}}\frac{C_{\nu}^{2}(r)}{\nu!}
\biggr)^{-1}\sum_{l\geq\kappa+1}\sum
_{\nu\in
N_{l}}\frac{C_\nu^2(r)}{\nu!}=\mathrm{o}\bigl(r^{\gamma/2}\bigr),\qquad r\to
\infty,
\end{equation}
where $\gamma\in(0,\min(\alpha,d-\alpha\kappa))$.

If there exists the limit distribution for at least one of the
random variables
\[
\frac{K_r}{\sqrt{ \operatorname{\mathbf{Var}} K_r}}\quad \mbox{and}\quad \frac
{K_{r,\kappa}}{\sqrt{\operatorname{ \mathbf{Var} } K_{r,\kappa}}},
\]
then the limit distribution of the other random variable exists too and
the limit distributions coincide when $r\to\infty$.
\end{theorem}

\begin{rem} If $G_r(w)\in\mathbf{L}_{2}(R^{p},\phi(\llVert  w\rrVert )  \,\mathrm{d}w)$ does
not depend on $r$ and has Hermitian rank $\kappa$, then (\ref{con1}) is
satisfied.
\end{rem}

\begin{rem}\label{varK} In many cases it is much easier to compute
$\operatorname{\mathbf{Var}}  K_{r,\kappa}$ than $\operatorname{\mathbf{Var}}  K_{r}$.
Using the property
$\lim_{r\to\infty}{\operatorname{\mathbf{Var}}  K_{r}}/{\operatorname{\mathbf{Var}}  K_{r,\kappa
}} = 1,
$
we can change the statement of Theorem~\ref{th4} as follows:
under the assumptions of Theorem~\ref{th4} limit distributions of the
random variables
${K_r}/{\sqrt{ \operatorname{\mathbf{Var}}  K_{r,\kappa}}}$ and ${K_{r,\kappa}}/{\sqrt { \operatorname{\mathbf{Var}}  K_{r,\kappa}}}$
coincide when $r\to\infty$.
\end{rem}
%
\begin{assumption}\label{ass2} $\eta_1(x)$ has a spectral density
$f(\llVert  \lambda\rrVert  )$, $\lambda\in\mathbb{R}^d$, such that
%
%
\begin{equation}
\label{f} f\bigl(\llVert \lambda\rrVert \bigr)\sim c_2(d,\alpha)
\llVert \lambda\rrVert ^{\alpha
-d}L \biggl( \frac1{\llVert \lambda\rrVert }
\biggr), \qquad\llVert \lambda \rrVert \to0,
\end{equation}
where $0<\alpha<d$ and
\[
c_2(d,\alpha):=\frac{\Gamma ( {(d-\alpha)}/2 ) }{2^\alpha
\pi^{d/2}\Gamma
( \alpha/2 ) }.
\]
\end{assumption}

\begin{rem}
If $f(\cdot)$ is decreasing in a neighbourhood of zero and continuous
for all $\lambda\neq0$, then by Tauberian Theorem~4 \cite{leoole} the
statement $\mathrm{B}(\llVert  x\rrVert )=\mathbf{E}\eta_1 (0)  \eta
_1(x)= \llVert  x\rrVert
^{-\alpha}L(\llVert  x\rrVert )$ implies Assumption~\ref{ass2}.
A much more detailed discussion of relations between Assumption~\ref
{ass1} and \ref{ass2} can be found in \cite{leoole,ole}.
\end{rem}

Note that then the field possesses the spectral representation
\[
\eta_1 (x)=\int_{\mathbb{R}^d}\mathrm{e}^{\mathrm{i}\langle\lambda,x\rangle}\sqrt{f
\bigl( \llVert \lambda\rrVert \bigr) }W(\mathrm{d}\lambda),
\]
where $W(\cdot)$ is the complex Gaussian white noise random measure on
$\mathbb{R}^d$.

Let
%
%
\begin{equation}
\label{ind} \mathcal{K}(x):=\int_{\Delta}\mathrm{e}^{\mathrm{i}\langle x,u\rangle} \,\mathrm{d}u, \qquad x
\in\mathbb{R}^d.
\end{equation}

\begin{theorem}\label{th5} Let $\eta_1(x)$, $x\in\mathbb{R}^d$, be a
homogeneous isotropic Gaussian random
field with $\mathbf{E}\eta_1(x)=0$. If Assumptions~\ref{ass1} and
\ref{ass2} hold, $\alpha\in(0,d/\kappa)$, and $\kappa\ge1$,
then for $r\to\infty$ the finite-dimensional distributions of
\[
X_{\kappa,r}:=r^{(\kappa\alpha)/2-d}L^{-\kappa/2}(r)\int_{\Delta
(r)}H_\kappa
\bigl(\eta_1 (x)\bigr) \,\mathrm{d}x
\]
converge weakly to the finite-dimensional distributions of
%
%
\begin{equation}
\label{Xk} X_\kappa:=c_2^{\kappa/2}(d,\alpha) \int
_{\mathbb{R}^{d\kappa
}}^{{\prime}}\mathcal{K}(\lambda_1+
\cdots +\lambda_\kappa) \frac{W(\mathrm{d}\lambda_1)\cdots W(\mathrm{d}\lambda_\kappa)}{\llVert  \lambda
_1\rrVert  ^{(d-\alpha)/2}\cdots\llVert  \lambda_\kappa\rrVert
^{(d-\alpha)/2}},
\end{equation}
where $\int_{R^{d\kappa}}^{{\prime}}$ denotes the multiple Wiener--It\^
{o} integral.
\end{theorem}

The following result shows that $X_\kappa$ is correctly defined and
${\mathbf E}X_\kappa^2<\infty$.
%
\begin{lemma}\label{finint}
If $\tau_1,\ldots,\tau_\kappa$, $\kappa\ge1$, are such positive
constants, that $\sum_{i=1}^\kappa\tau_i <d$, then
%
%
\begin{equation}
\label{finv} \int_{\mathbb{R}^{d\kappa}}\bigl|\mathcal{K}(\lambda_1+
\cdots +\lambda_\kappa)\bigr|^2 \frac{\mathrm{d}\lambda_1\cdots \mathrm{d}\lambda_\kappa}{\llVert
\lambda
_1\rrVert  ^{d-\tau_1}\cdots\llVert  \lambda_\kappa\rrVert  ^{d-\tau
_\kappa}}<\infty.
\end{equation}

If $\tau_1=\cdots=\tau_\kappa=\alpha$, $\alpha\in(0,d/\kappa)$, then we
will use the following notation
\[
c_3(\kappa,d,\alpha):=\int_{\mathbb{R}^{d\kappa}}\bigl|\mathcal{K}(
\lambda _1+\cdots +\lambda_\kappa)\bigr|^2
\frac{\mathrm{d}\lambda_1\cdots \mathrm{d}\lambda_\kappa}{\llVert
\lambda
_1\rrVert  ^{d-\alpha}\cdots\llVert  \lambda_\kappa\rrVert
^{d-\alpha}} .
\]
\end{lemma}

\begin{rem}
It is not difficult to adapt Theorem~\ref{th5} for the case of
stochastic processes and obtain self-similar limit processes, consults
\cite{iv,leo1,leoole0,mes}.
For $\kappa=2$, the limit random variable $X_2$ in Theorem~\ref{th5}
plays an analogous role to the Rosenblatt distribution, see \cite{ta0}.
\end{rem}

\begin{example}
If $\Delta$ is the ball $v(1)$, then
\[
\mathcal{K}(x)=\int_{v(1)}\mathrm{e}^{\mathrm{i}\langle x,u\rangle } \,\mathrm{d}u=(2
\pi)^{d/2} \frac{J_{d/2}(\|x\|)}{\llVert x\rrVert ^{d/2}},\qquad x\in\mathbb{R}^d,
\]
and we obtain the result from \cite{iv}, Section~2.10, with $t=1$, that is,
\[
X_\kappa=(2\pi)^{d/2}c_2^{\kappa/2}(d,\alpha)
\int_{\mathbb{R}^{d\kappa
}}^{{\prime}}\frac{J_{d/2}(\llVert  \lambda_1+\cdots
+\lambda_\kappa\rrVert  )}{\llVert  \lambda_1+\cdots+\lambda_\kappa
\rrVert
^{d/2}}
\frac{W(\mathrm{d}\lambda_1)\cdots W(\mathrm{d}\lambda_\kappa)}{\llVert  \lambda
_1\rrVert  ^{(d-\alpha)/2}\cdots\llVert  \lambda_\kappa\rrVert
^{(d-\alpha)/2}%
}.
\]
\end{example}

\begin{example} Let us consider $\eta(x)$ with uncorrelated
identically distributed components possessing covariance functions of
the form
\[
B_{jj} \bigl( \llVert x\rrVert \bigr) = \bigl( 1+\llVert x\rrVert
^{{}\sigma} \bigr) ^{-\theta},\qquad  \sigma\in ( 0,2 ] , \theta>0, j=1,\ldots,p.
\]
The above is known as the generalized Linnik covariance function.
Cauchy field in the simulation results of Section~\ref{sec9} is an
important particular case of this model.

If $\sigma\theta\kappa> d$, $\kappa\ge1$, then $\eta(x)$ is a
weakly dependent random field which satisfies the assumptions of
Section~\ref{sec4}, that is, $\psi(x)=B_{11} ( \llVert  x\rrVert
)\in\mathbf{L}_\kappa(\mathbb{R}^d)$ and $\psi(x)\le1$ for all
$x\in\mathbb{R}^d$. If $\sigma\theta< d$, then we have the strongly
dependent case and Assumptions~\ref{ass1} and \ref{ass2}
hold, see \cite{leoole} and references therein.
\end{example}

In the next two theorems, we apply the general results to study the
sojourn measure of strongly dependent Fisher--Snedecor and Student
random fields above a constant level, that is, $a(r)\equiv a$.
The following theorem demonstrates that for Student random fields, even
in the case of strong dependence, we have a normal limit law. However,
for the strong\-ly dependent case the normalization is different from
$r^{-d/2}$ in Theorem~\ref{thS}.

\begin{theorem}\label{th6} Let $\eta(x)=[\eta_{1}(x),\ldots,\eta
_{n+1}(x)]^{\prime}$, $x\in\mathbb{R}^{d}$, satisfy Assumption~\ref{ass1} for $\alpha\in(0,d)$, and Assumption~\ref{ass2} hold
for the spectral density of each component $\eta_{j}(\cdot)$. Then the
random variable
\[
U_r(n,\alpha):=\sqrt{2\pi} \bigl(1+a^2/n
\bigr)^{n/2}\frac{M_{r} \{
T_{n} \}-\llvert  \Delta\rrvert  r^{d} (1/2-1/2
(1-I_{{n}/{(n+a^2)}} ({n}/{2},{1}/{2} ) )\cdot
\operatorname{sgn}(a) )}{r^{d-\alpha/2}L^{1/2}(r)\sqrt{c_2(d,\alpha)
c_3(1,d,\alpha)}}
\]
is asymptotically $\mathcal{N}(0,1)$, as $r\to\infty$.
\end{theorem}

Contrary to the Student case, for strongly dependent Fisher--Snedecor
random fields we obtain a non-normal limit law.

\begin{theorem}\label{th7} Let $\eta(x)=[\eta_{1}(x),\ldots,\eta
_{n+m}(x)]^{\prime}$, $x\in\mathbb{R}^{d}$, satisfy Assumption~\ref{ass1} for $\alpha\in(0,d/2)$, and Assumption~\ref{ass2}
hold for the spectral density of each component $\eta_{j}(\cdot)$.
Then, for $r\to\infty$,
the distribution of the random variable
\[
U_r(m,n,\alpha):=\frac{M_{r} \{ F_{m,n} \}-\llvert  \Delta\rrvert  r^{d} (1-I_{{ma}/{(n+ma)}} ({m}/{2},{n}/{2}
) )}{c_4(a,n,m) r^{d-\alpha}L(r)}
\]
converges to the distribution of the random variable
\[
R(m,n):=\frac{X_{2,1}+\cdots+X_{2,m}}{m}-\frac{X_{2,m+1}+\cdots+X_{2,m+n}}{n},
\]
where $X_{2,j}$, $j=1,\ldots,m+n$, are independent copies of the random
variable $X_{2}$ defined by~(\ref{Xk}),
\[
c_4(a,n,m):= \frac{(ma/n)^{m/2} \Gamma ( (m+n)/2 )
}{(1+ma/n)^{(m+n)/2} \Gamma(n/2) \Gamma (m/2 )}.
\]
\end{theorem}

Now we generalize the previous results to the increasing level $a(r)\to
+\infty$, as $r\to+\infty$.

\begin{theorem}\label{th8} Let $\eta(x)=[\eta_{1}(x),\ldots,\eta
_{n+1}(x)]^{\prime}$, $x\in\mathbb{R}^{d}$, satisfy Assumption~\ref{ass1} for $\alpha\in(0,d)$, and Assumption~\ref{ass2} hold
for the spectral density of each component $\eta_{j}(\cdot)$. If
$a(r)=\mathrm{o}(r^{\gamma/2n})$, $\gamma\in(0,\min(\alpha,d-\alpha))$, $r\to
\infty$, then
the random variable
\[
\sqrt{2\pi} \bigl(1+a(r)^2/n \bigr)^{n/2}\frac{M_{r} \{ T_{n} \}
-\llvert  \Delta\rrvert  r^{d}I_{{n}/{(n+a^2(r))}} (
{n}/{2},{1}/{2} )}{r^{d-\alpha/2}L^{1/2}(r)\sqrt{c_2(d,\alpha
) c_3(1,d,\alpha)}}
\]
is asymptotically $\mathcal{N}(0,1)$.
\end{theorem}

\begin{theorem}\label{th9} Let $\eta(x)=[\eta_{1}(x),\ldots,\eta
_{n+m}(x)]^{\prime}$, $x\in\mathbb{R}^{d}$, satisfy Assumption~\ref{ass1} for $\alpha\in(0,d/2)$, and Assumption~\ref{ass2}
hold for the spectral density of each component $\eta_{j}(\cdot)$. If
$a(r)=\mathrm{o}(r^{\gamma/n})$, $\gamma\in(0,\min(\alpha,d-\alpha))$, $r\to
\infty$, then the distribution of the random variable
\[
\frac{M_{r} \{ F_{m,n} \}-\llvert  \Delta\rrvert  r^{d}
(1-I_{{ma(r)}/{(n+ma(r))}} ({m}/{2},{n}/{2} )
)}{c_4(a(r),n,m) r^{d-\alpha}L(r)}
\]
converges to the distribution of the random variable $R(m,n)$
defined in Theorem~\ref{th7}.
\end{theorem}

The following theorems illustrate how to extend the obtained results to
long range dependent vector fields which components may be
cross-correlated, consult the pioneering papers \cite{mae1,mae2,ta1}
on similar vector Gaussian process results. Such cross-correlated
random fields may be useful in positron emission tomography studies to
identify brain activated regions. In many cases, the activation is so
small that the experiment must be repeated several times and the scan
results are averaged to improve the signal-to-noise ratio. The
cross-correlated components $\eta_j(x)$, $j=1,\ldots,p$, can be
interpreted as repeated imaged slices in scans of the same subject. If
the stationarity assumption is in doubt, Student and Fisher--Snedecor
random fields were proposed to test regional changes, consult \cite{cao1,wor}.

We use the previous notation $M_r\{T_n\}$ and $M_{r}\{ F_{m,n}\}$,
but replace independent components of $\eta(\cdot)$ in the
definitions~\ref{def3} and \ref{def4} by components of cross-correlated
random fields. Note, that the functional $M_r\{T_n\}$ ($M_{r}\{
F_{m,n}\}$) takes the same value on the class of fields $\{C\eta(x), C>0\}$. Therefore, we study only the cases where $\det(\mathbf{E}\eta
(0)  \eta(0)^{\prime})=1$.

\begin{assumption}\label{ass3}
Let $\eta(x)=[\eta_{1}(x),\ldots,\eta_{p}(x)]^{\prime}$, $x\in
\mathbb{R}^{d}$, be a vector homogeneous isotro\-pic zero mean Gaussian
random field such that
\[
\mathbf{B}\bigl(\llVert x\rrVert \bigr)=\mathbf{E}\eta(0) \eta(x)^{\prime
}=
\mathcal{A}\cdot\llVert x\rrVert ^{-\alpha}L\bigl(\llVert x\rrVert \bigr),\qquad
\alpha\in(0,d/\kappa), \kappa\ge1,
\]
where $\mathcal{A}$ is a $p\times p$ positive-semidefinite symmetric
orthogonal matrix, and
Assumption~\ref{ass2} hold for the spectral density of each
component of the field $\tilde{\eta}:=\mathcal{A}^{-1/2}\eta$.
\end{assumption}

Note that, by the definition of $\mathcal{A}$, there exists the square
root of $\mathcal{A}^{-1}$, that is, the positive-semidefinite
orthogonal matrix $\mathcal{A}^{-1/2}$, such that $\mathcal
{A}^{-1/2}\mathcal{A}^{-1/2}=\mathcal{A}^{-1}$. In what follows, we
denote $\mathcal{A}^{-1/2}:=(a_{ij})_{1\leq i,j\leq p}$.

\begin{theorem}\label{th10} If $\eta(x)=[\eta_{1}(x),\ldots,\eta
_{n+1}(x)]^{\prime}$, $x\in\mathbb{R}^{d}$, satisfies Assumption~\ref{ass3} for $\kappa=1$,
then $U_r(n,\alpha)$ defined in Theorem~\ref{th6}
is asymptotically $\mathcal{N}(0,1)$, as $r\to\infty$.
\end{theorem}

For the Fisher--Snedecor random field, we only consider the case of a
block diagonal matrix $\mathcal{A}$. It is also possible to derive
similar results for arbitrary $\mathcal{A}$, but for such cases we need
a generalization of Theorem~\ref{th5} about the asymptotic behaviour of
the bivariate functionals $\int_{\Delta(r)}\eta_j (x)\eta_l (x) \,\mathrm{d}x$
(consult \cite{ta1} for $d=1$), which is beyond the scope of this paper.
%
\begin{theorem}\label{th11} Let $\eta(x)=[\eta_{1}(x),\ldots,\eta
_{n+m}(x)]^{\prime}$, $x\in\mathbb{R}^{d}$, satisfy Assumption~\ref{ass3} for $\kappa=2$
and $\mathcal{A}= \bigl[
{\mathcal{A}_1 \atop 0} \enskip
{ 0 \atop \mathcal{A}_2}
\bigr]$, where $\mathcal{A}_1$ and $\mathcal{A}_2$ are $m\times m$ and
$n\times n$ matrices, respectively.
Then, for $r\to\infty$,
the distribution of the random variable $U_r(m,n,\alpha)$
converges to the distribution of the random variable $R(m,n)$,
where $U_r(m,n,\alpha)$ and $R(m,n)$ are defined in Theorem~\ref{th7}.
\end{theorem}

\section{Proofs of the results of Section~\texorpdfstring{\protect\ref{sec6}}{7}}\label{proofs2}

\begin{pf*}{Proof of Theorem~\ref{th4}}Let
\[
V_r:=\sum_{l\geq\kappa+1} \sum
_{\nu\in N_l}\frac{C_\nu(r)}{\nu!} \int_{\Delta(r)}e_\nu
\bigl(\eta(x)\bigr)\,\mathrm{d}x,
\]
then by Lemma~\ref{lem1}
\[
\operatorname{\mathbf{Var}} K_r= \operatorname{\mathbf{Var}} K_{r,\kappa}+ \operatorname{\mathbf{Var}}
V_r.
\]
By (\ref{varg}) and (\ref{dint})
\begin{eqnarray*}
\operatorname{\mathbf{Var}} K_{r,\kappa}&=& \sum_{\nu\in N_\kappa}
\frac{C_\nu
^2(r)}{\nu!} \int_{\Delta(r)}\int_{\Delta(r)}
\llVert x-y\rrVert ^{-\alpha\kappa} L^\kappa \bigl(\llVert x-y\rrVert
\bigr) \,\mathrm{d}x \,\mathrm{d}y
\\
&=&|\Delta|^2r^{2 d-\alpha\kappa}\sum_{\nu\in N_\kappa}
\frac{C_\nu
^2(r)}{\nu!} \int_0^{\mathrm{diam} \{ \Delta
\}} z^{-\alpha\kappa}
L^\kappa (rz )\psi_{\Delta}(z)\,\mathrm{d}z.
\end{eqnarray*}

If $\alpha\in(0,d/\kappa)$, then by asymptotic properties of
integrals of slowly varying functions (see Theorem~2.7 \cite{sen}) we get
\begin{eqnarray*}
\operatorname{\mathbf{Var}} K_{r,\kappa}&=&c_1(\kappa,\alpha,\Delta) |
\Delta|^2 \sum_{\nu\in N_\kappa}\frac{C_\nu^2(r)}{\nu!}
r^{2d-\kappa\alpha}L^\kappa(r) \bigl(1+\mathrm{o}(1)\bigr),\qquad r\to\infty,
\end{eqnarray*}
where
\[
c_1(\kappa,\alpha,\Delta):=\int_0^{\mathrm{diam} \{ \Delta
\}}
z^{-\alpha\kappa}\psi_{\Delta}(z)\,\mathrm{d}z.
\]

Similar to $\operatorname{\mathbf{Var}}  K_{r,\kappa}$ we obtain
\[
\operatorname{\mathbf{Var}} V_r = |\Delta|^2r^{2d}\sum
_{l\geq\kappa+1}\sum_{\nu\in
N_{l}}
\frac{C_\nu^2(r)}{\nu!} \int_0^{r\cdot \mathrm{diam} \{ \Delta
\}} z^{-\alpha l}
L^l (z )\psi_{\Delta(r)}(z)\,\mathrm{d}z.
\]
It follows from $z^{-\alpha} L (z )\in[0,1]$, $z\ge0$, that
\begin{eqnarray*}
\operatorname{\mathbf{Var}} V_r &\leq& |\Delta|^2r^{2d-(\kappa+1)\alpha}\sum
_{l\geq
\kappa+1}\sum_{\nu\in
N_{l}}
\frac{C_\nu^2(r)}{\nu!} \int_0^{\mathrm{diam} \{ \Delta
\}} z^{-\alpha(\kappa+1)}
L^{\kappa+1} (rz )\psi _{\Delta}(z)\,\mathrm{d}z
\\
&=& |\Delta|^2r^{2d-\kappa\alpha}L^{\kappa}(r)\sum
_{l\geq\kappa+1}\sum_{\nu\in
N_{l}}
\frac{C_\nu^2(r)}{\nu!} \int_0^{\mathrm{diam} \{ \Delta
\}} z^{-\alpha\kappa}
\frac{L^{\kappa} (rz )}{L^{\kappa
}(r)} \frac{L (rz )}{(rz)^{\alpha}}\psi_{\Delta}(z)\,\mathrm{d}z.
\end{eqnarray*}

Let us split the above integral into two parts $I_1$ and $I_2$ with the
ranges of integration $[0,r^{-\beta}]$ and $(r^{-\beta},\mathrm{diam} \{
\Delta
\}]$, respectively, where $\beta\in(0,1)$.

As $z^{-\alpha} L (z )\in[0,1]$, $z\ge0$, we can estimate the
first integral as follows
\begin{eqnarray}\label{int}
I_1&\le&\int_0^{r^{-\beta}}
z^{-\alpha\kappa} \frac{L^{\kappa}
(rz )}{L^{\kappa}(r)} \psi_{\Delta}(z)\,\mathrm{d}z\le
\frac{\sup_{0\le s\le r^{1-\beta}}s^{\delta}L^{\kappa} (s
)}{r^{\delta}L^{\kappa}(r)} \int_0^{r^{-\beta}}
z^{-\delta}z^{-\alpha
\kappa} \psi_{\Delta}(z)\,\mathrm{d}z
\nonumber
\\[-8pt]
\\[-8pt]
\nonumber
 &\le& \biggl(\frac{\sup_{0\le s\le r}s^{\delta/k}L (s )}{r^{\delta
/k}L(r)} \biggr)^{\kappa} \int
_0^{r^{-\beta}} z^{-\delta}z^{-\alpha\kappa
}
\psi_{\Delta}(z)\,\mathrm{d}z.
\end{eqnarray}

By Theorem~1.5.3 \cite{bin} and the definition of slowly varying functions
\[
\lim_{r\to\infty}\frac{\sup_{0\le s\le r}s^{\delta/k}L (s
)}{r^{\delta/k}L(r)}=1.
\]

By (\ref{dint}), we can estimate the integral in (\ref{int}) as follows
\begin{eqnarray}\label{I1}
\int_0^{r^{-\beta}} z^{-\delta}z^{-\alpha\kappa}
\psi_{\Delta
}(z)\,\mathrm{d}z&=&\llvert \Delta\rrvert ^{-2}\int
_{\Delta}\int_{\Delta}\chi\bigl(\llVert x-y
\rrVert \le r^{-\beta}\bigr)\llVert x-y\rrVert ^{-(\delta+\alpha\kappa
)} \,\mathrm{d}x \,\mathrm{d}y
\nonumber
\\[-8pt]
\\[-8pt]
\nonumber
&\le&\llvert \Delta\rrvert ^{-1}\int_{0}^{r^{-\beta}}
\rho^{d-(1+\delta+\alpha\kappa)} \,\mathrm{d}\rho=\frac{r^{-\beta
(d-(\delta+\alpha\kappa))}}{(d-(\delta+\alpha\kappa))\llvert \Delta\rrvert }.
\end{eqnarray}

For the second integral, we obtain
\[
I_2\le \frac{\sup_{r^{1-\beta}\le s\le r\cdot \mathrm{diam} \{ \Delta
\}}s^{\delta}L^{\kappa} (s )}{r^{\delta}L^{\kappa
}(r)}\cdot\sup_{r^{1-\beta}\le s\le r\cdot \mathrm{diam} \{ \Delta
\}}
\frac{L (s )}{s^\alpha}\cdot\int_0^{\mathrm{diam} \{
\Delta \}}
z^{-(\delta+\alpha\kappa)} \psi_{\Delta}(z)\,\mathrm{d}z.
\]

Using Theorem~1.5.3 \cite{bin}, we conclude that
\begin{eqnarray*}
\lim_{r\to\infty}\frac{\sup_{r^{1-\beta}\le s\le r\cdot \mathrm{diam} \{
\Delta
\}}s^{\delta}L^{\kappa} (s )}{r^{\delta}L^{\kappa}(r)}&\le &\lim_{r\to\infty}
\frac{\sup_{0\le s\le r\cdot \mathrm{diam} \{ \Delta
\}}s^{\delta}L^{\kappa} (s )}{(r\cdot \mathrm{diam} \{ \Delta
\})^{\delta}L^{\kappa}(r\cdot \mathrm{diam} \{ \Delta \})}
\\
&&{}\times\lim_{r\to\infty}\frac{\mathrm{diam}^{\delta} \{ \Delta
\}L^{\kappa}(r\cdot \mathrm{diam} \{ \Delta \})}{L^{\kappa
}(r)}=\mathrm{diam}^{\delta} \{
\Delta \}.
\end{eqnarray*}

By Proposition~1.3.6 and Theorem~1.5.3 \cite{bin}, it follows that
%
%
\begin{equation}
\label{I2}\sup_{r^{1-\beta}\le s\le r\cdot \mathrm{diam} \{
\Delta \}}\frac{L (s )}{s^\alpha}\le\frac{\sup_{s\ge
r^{1-\beta}}s^{-\alpha}L (s )}{r^{-\alpha(1-\beta)}L
(r^{1-\beta} )}
\cdot\frac{L (r^{1-\beta} )}{r^{\delta
(1-\beta)}} \cdot r^{(\delta-\alpha)(1-\beta)}=\mathrm{o}\bigl(r^{(\delta-\alpha
)(1-\beta)}\bigr).
\end{equation}

We can choose $\beta=1/2$ and make $\delta$ arbitrary close to 0. Then
by (\ref{I1}), (\ref{I2}), and condition (\ref{con1}) we obtain
\[
\lim_{r\to\infty}\frac{\operatorname{\mathbf{Var}}  V_r}{\operatorname{\mathbf{Var}}
K_{r}}=0 \quad\mbox{and}\quad \lim
_{r\to\infty}\frac{\operatorname{\mathbf{Var}}  K_{r}}{\operatorname{\mathbf{Var}}
K_{r,\kappa}} = 1.
\]

Thus,
\[
\lim_{r\to\infty} \mathbf{E} \biggl(\frac{K_{r}}{\sqrt{ \operatorname{\mathbf{Var}}
K_{r}}}-
\frac{K_{r,\kappa}}{\sqrt{\operatorname{ \mathbf{Var}}  K_{r,\kappa}}} \biggr)^2=\lim_{r\to\infty}
\frac{\mathbf{E} (V_r+ (1-\sqrt{
{\operatorname{\mathbf{Var}} K_{r}}/{\operatorname{\mathbf{Var}} K_{r,\kappa}}} )K_{r,\kappa
} )^2}{\operatorname{\mathbf{Var}} K_{r}} =0,
\]
which completes the proof.
\end{pf*}

\begin{pf*}{Proof of Lemma~\ref{finint}}
Definition (\ref{ind}) yields $\mathcal{K}(\cdot)\in\mathbf{L}_\infty(
\mathbb{R}^d)$ and by the Plancherel theorem $\mathcal{K}(\cdot)\in
\mathbf{L}_2(\mathbb{R}^d)$.
Hence, the statement of the lemma is valid for $\kappa= 1$.
For $\kappa> 1$, we can obtain (\ref{finv}) by the recursive estimation
routine and the change of variables $\tilde{\lambda}_{\kappa-1}={\lambda
_{\kappa-1}}/{\llVert  u\rrVert }$:
\begin{eqnarray*}
&&\int_{\mathbb{R}^{d\kappa}}\bigl|\mathcal{K}(\lambda_1+\cdots +
\lambda_{\kappa})\bigr|^2 \frac{\mathrm{d}\lambda_1\cdots \mathrm{d}\lambda_\kappa}{\llVert  \lambda
_1\rrVert  ^{d-\tau_1}\cdots\llVert  \lambda_\kappa\rrVert  ^{d-\tau
_\kappa}}\\
&&\quad= \int
_{\mathbb{R}^{d(\kappa-1)}}\bigl|\mathcal{K}(\lambda_1+\cdots +
\lambda_{\kappa-2}+u)\bigr|^2\\
&&\qquad{}\times\int_{\mathbb{R}^{d}} \frac{d\lambda_{\kappa-1}}{\llVert
\lambda_{\kappa-1}\rrVert  ^{d-\tau_{\kappa-1} }\llVert u- \lambda
_{\kappa-1}\rrVert  ^{d-\tau_{\kappa}}}\cdot\frac{\mathrm{d}\lambda_1\cdots
\mathrm{d}\lambda_{\kappa-2}  \,\mathrm{d}u}{\llVert  \lambda
_1\rrVert  ^{d-\tau_{1} }\cdots\llVert  \lambda_{\kappa-2}\rrVert
^{d-\tau_{\kappa-2} }}
\\
&&\quad= \int_{\mathbb{R}^{d(\kappa-1)}}\frac{|\mathcal{K}(\lambda_1+\cdots
+\lambda_{\kappa-2}+u)|^2 \,\mathrm{d}\lambda_1\cdots \mathrm{d}\lambda_{\kappa-2}}{\llVert  \lambda
_1\rrVert  ^{d-\tau_{1} }\cdots\llVert  \lambda_{\kappa-2}\rrVert
^{d-\tau_{\kappa-2} }\llVert  u\rrVert ^{d-\tau_{\kappa-1}-\tau_{\kappa}
}}\\
&&\qquad{}\times \int_{\mathbb{R}^{d}}
\frac{\mathrm{d}\tilde{\lambda}_{\kappa-1}}{\llVert
\tilde{\lambda}_{\kappa-1}\rrVert  ^{d-\tau_{\kappa-1} }\llVert
{u}/{\llVert  u\rrVert }- \tilde{\lambda}
_{\kappa-1}\rrVert  ^{d-\tau_{\kappa}}} \,\mathrm{d}u
\\
&&\quad
\le C \int_{\mathbb{R}^{d(\kappa-1)}}\bigl|\mathcal{K}(\lambda_1+\cdots
+\lambda_{\kappa-2}+u)\bigr|^2 \frac{\mathrm{d}\lambda_1\cdots \mathrm{d}\lambda_{\kappa
-2}  \,\mathrm{d}u}{\llVert  \lambda
_1\rrVert  ^{d-\tau_{1} }\cdots\llVert  \lambda_{\kappa-2}\rrVert
^{d-\tau_{\kappa-2} }\llVert  u\rrVert ^{d-\tau_{\kappa-1}-\tau_{\kappa} }}
\\
&&\quad\le\cdots\\
&&\quad\le C \int_{\mathbb{R}^{d}}\bigl|\mathcal{K}(u)\bigr|^2
\frac{\mathrm{d}u}{\llVert
u\rrVert ^{d-\sum_{i=1}^\kappa\tau_i}}<\infty.
\end{eqnarray*}
\upqed\end{pf*}

\begin{pf*}{Proof of Theorem~\ref{th5}}
Using the self-similarity of Gaussian white noise, namely $W (
C \,\mathrm{d}\lambda ) \stackrel{\mathcal{D}}{=}C^{d/2}W ( \mathrm{d}\lambda
)$, and the It\'o formula \cite{dob}
\[
H_\kappa\bigl(\eta_1 (x)\bigr)=\int_{\mathbb{R}^{d\kappa}}^{\prime
}\mathrm{e}^{\mathrm{i}\langle\lambda_1+\cdots+\lambda
_\kappa,x\rangle}
\Biggl\{ \prod_{j=1}^\kappa\sqrt{f(
\lambda_j)} \Biggr\} W(\mathrm{d}\lambda _1)\cdots W(\mathrm{d}
\lambda_\kappa)
\]
we obtain
\begin{eqnarray*}
&&X_{\kappa,r}\stackrel{\mathcal{D}} {=} c_2^{\kappa/2}(d,
\alpha)\int_{\mathbb{R}^{d\kappa}}^{{\prime}}\mathcal {K}(
\lambda_1+\cdots+\lambda_\kappa)Q_r(
\lambda_1,\ldots,\lambda _\kappa)\frac{W(\mathrm{d}\lambda
_1)\cdots W(\mathrm{d}\lambda_\kappa)}{\llVert  \lambda_1\rrVert  ^{(d-\alpha
)/2}\cdots
\llVert  \lambda_\kappa\rrVert  ^{(d-\alpha)/2}},
\end{eqnarray*}
where
\[
Q_r(\lambda_1,\ldots,\lambda_\kappa):
=r^{\kappa(\alpha
-d)/2}L^{-\kappa/2}(r) c_2^{-\kappa/2}(d,\alpha)
\Biggl[ \prod_{j=1}^\kappa\llVert
\lambda_j\rrVert ^{d-\alpha}f \biggl( \frac{\llVert  \lambda_j\rrVert
}r \biggr)
\Biggr] ^{1/2}.
\]

By the isometry property of multiple stochastic integrals
\[
R_r:=\frac{\mathbb{E}\llvert  X_{\kappa,r}-X_\kappa\rrvert ^2}{c_2^{\kappa
}(d,\alpha)}=\int_{\mathbb{R}^{d\kappa}}
\frac{|\mathcal{K}(\lambda
_1+\cdots+\lambda_\kappa)|^2  (Q_r(\lambda_1,\ldots,\lambda
_\kappa)-1 )^2}{\llVert  \lambda_1\rrVert  ^{d-\alpha}\cdots
\llVert  \lambda_\kappa\rrVert  ^{d-\alpha}} \,\mathrm{d}\lambda _1\cdots \mathrm{d}\lambda_\kappa.
\]

Using (\ref{f}) and properties of slowly varying functions we conclude
that $Q_r(\lambda_1,\ldots,\lambda
_\kappa)$ converges pointwise to 1, when $r\to\infty$.
Hence, by Lebesgue's dominated convergence theorem the integral
converges to zero if there is some integrable function which dominates
integrands for all~$r$.

Let us split $\mathbb{R}^{d\kappa}$ into the regions
\[
B_\mu:=\bigl\{(\lambda_1,\ldots,\lambda_\kappa)
\in\mathbb{R}^{d\kappa}\dvt \Vert \lambda_j\Vert \le1,
\mbox{ if } \mu_j=-1, \mbox{ and } \Vert \lambda_j\Vert >
1, \mbox{ if } \mu_j=1, j=1,\ldots,\kappa\bigr\},
\]
where $\mu=(\mu_1,\ldots,\mu_\kappa)\in\{-1,1\}^\kappa$ is a binary
vector of length $\kappa$.
Then we can represent the integral $R_r$ as
\[
R_r:=\bigcup_{\mu\in\{-1,1\}^\kappa}\int
_{B_\mu}\bigl|\mathcal{K}(\lambda _1+\cdots+
\lambda_\kappa)\bigr|^2 \bigl(Q_r(
\lambda_1,\ldots,\lambda _\kappa)-1 \bigr)^2
\frac{\mathrm{d}\lambda
_1\cdots\mathrm{ d}\lambda_\kappa}{\llVert  \lambda_1\rrVert  ^{d-\alpha}\cdots
\llVert  \lambda_\kappa\rrVert  ^{d-\alpha}}.
\]

If $(\lambda_1,\ldots,\lambda_\kappa)\in B_\mu$ we estimate the
integrand as follows
\begin{eqnarray*}
&&\frac{|\mathcal{K}(\lambda_1+\cdots+\lambda_\kappa)|^2
(Q_r(\lambda_1,\ldots,\lambda_\kappa)-1 )^2}{\llVert  \lambda
_1\rrVert  ^{d-\alpha}\cdots
\llVert  \lambda_\kappa\rrVert  ^{d-\alpha}}\\
&&\quad\le\frac{2  |\mathcal
{K}(\lambda_1+\cdots+\lambda_\kappa)|^2 }{\llVert  \lambda_1\rrVert
^{d-\alpha}\cdots
\llVert  \lambda_\kappa\rrVert  ^{d-\alpha}} \bigl(Q^2_r(
\lambda _1,\ldots,\lambda_\kappa)+1 \bigr)
\\
&&\quad= \frac{2  |\mathcal{K}(\lambda_1+\cdots+\lambda_\kappa)|^2 }{\llVert  \lambda_1\rrVert  ^{d-\alpha}\cdots
\llVert  \lambda_\kappa\rrVert  ^{d-\alpha}} \Biggl(\prod_{j=1}^\kappa
\Vert \lambda_j\Vert ^{\mu_j\delta}\cdot\prod
_{j=1}^\kappa\frac{ ({r}/{\Vert \lambda_j\Vert } )^{\mu_j\delta
}L ({r}/{\Vert \lambda_j\Vert } )}{r^{\mu_j\delta}L(r)} +1 \Biggr)
\\
&&\quad\le\frac{2 |\mathcal{K}(\lambda_1+\cdots+\lambda_\kappa)|^2 }{\llVert  \lambda_1\rrVert  ^{d-\alpha}\cdots
\llVert  \lambda_\kappa\rrVert  ^{d-\alpha}} \Biggl(1+ \prod_{j=1}^\kappa
\llVert \lambda_1\rrVert ^{\mu_j\delta}\cdot \sup
_{(\lambda_1,\ldots,\lambda_\kappa)\in B_\mu}\prod_{j=1}^\kappa
\frac{ ({r}/{\Vert \lambda_j\Vert } )^{\mu_j\delta
}L ({r}/{\Vert \lambda_j\Vert } )}{r^{\mu_j\delta}L(r)} \Biggr),
\end{eqnarray*}
where $\delta$ is an arbitrary positive number.
By Theorem~1.5.3 \cite{bin}
\begin{eqnarray*}
\lim_{r\to\infty}\frac{\sup_{\Vert \lambda_j\Vert \le1} (
{r}/{\Vert \lambda_j\Vert } )^{-\delta}L ({r}/{\Vert \lambda
_j\Vert } )}{r^{-\delta}L(r)}&=&\lim_{r\to\infty}
\frac{\sup_{z\ge
r}z^{-\delta}L (z )}{r^{-\delta}L(r)}=1;
\\
\lim_{r\to\infty}\frac{\sup_{\Vert \lambda_j\Vert > 1} (
{r}/{\Vert \lambda_j\Vert } )^{\delta}L ({r}/{\Vert \lambda_j\Vert }
)}{r^{\delta}L(r)}&=&\lim_{r\to\infty}
\frac{\sup_{z\in[0,r]}z^{\delta
}L (z )}{r^{\delta}L(r)}=1.
\end{eqnarray*}

Therefore, there exists $r_0>0$ such that for all $r\ge r_0$ and
$(\lambda_1,\ldots,\lambda_\kappa)\in B_\mu$
\begin{eqnarray}\label{upper}
&&\frac{|\mathcal{K}(\lambda_1+\cdots+\lambda_\kappa)|^2
(Q_r(\lambda_1,\ldots,\lambda_\kappa)-1 )^2}{\llVert  \lambda
_1\rrVert  ^{d-\alpha}\cdots
\llVert  \lambda_\kappa\rrVert  ^{d-\alpha}}
\nonumber
\\
&&\quad\le\frac{2  |\mathcal
{K}(\lambda_1+\cdots+\lambda_\kappa)|^2 }{\llVert  \lambda_1\rrVert
^{d-\alpha}\cdots
\llVert  \lambda_\kappa\rrVert  ^{d-\alpha}}
 \\
 &&\qquad{}+2 C \frac{|\mathcal{K}(\lambda_1+\cdots+\lambda_\kappa)|^2 }{\llVert  \lambda_1\rrVert  ^{d-\alpha-\mu_1\delta}\cdots
\llVert  \lambda_\kappa\rrVert  ^{d-\alpha-\mu_\kappa\delta}}.\nonumber
\end{eqnarray}

By Lemma~\ref{finint}, if we chose $\delta\in (0,\min (\alpha
,{d}/{\kappa}-\alpha ) )$, the upper bound in (\ref{upper})
is an integrable function on each $B_\mu$ and hence on $\mathbb
{R}^{d\kappa}$ too.
By Lebesgue's dominated convergence theorem $\lim_{r\to\infty} \mathbf
{E}\llvert  X_{\kappa,r}-X_\kappa\rrvert ^2=0$, which completes the proof.
\end{pf*}

\begin{pf*}{Proof of Theorem~\ref{th6}}
For the function $\tilde{G} ( \cdot )$ given by (\ref{tG})
coefficients $C_{\nu}= 0$ for $\nu\in N_{1}\setminus\{(1,0,\ldots,0)\}
$. $C_{(1,0,\ldots,0)}$ is given by the formula
\begin{eqnarray}\label{C10}
C_{(1,0,\ldots,0)}&=&\int_{\mathbb{R}^{n+1}}\tilde {G}_{r}(w)e_{(1,0,\ldots,0)}(w)
\phi\bigl(\llVert w\rrVert \bigr)\,\mathrm{d}w\nonumber\\
&=&\int_{\mathbb{R}^{n+1}}\chi \biggl(
\frac{w_{1}}{\sqrt{{1}/{n}
(w_{2}^{2}+\cdots
+w_{n+1}^{2} )}}>a \biggr) w_{1}\prod_{j=1}^{n+1}
\frac{\mathrm{e}^{-{w_{j}^{2}}/{2}}}{\sqrt{2\pi
}}\,\mathrm{d}w_j
\nonumber
\\[-8pt]
\\[-8pt]
\nonumber
&=&\frac{2 \pi^{n/2}}{(2\pi)^{(n+1)/2}\Gamma(n/2)}\int_0^\infty
\rho ^{n-1} \mathrm{e}^{-{\rho^{2}}/{2}}\int_{|a|\rho/\sqrt{n}}^\infty
w_{1} \mathrm{e}^{-{w_{1}^{2}}/{2}} \,\mathrm{d}w_1 \,\mathrm{d}\rho
\\
&=&\frac{1}{\sqrt{2\pi} (1+a^2/n )^{n/2}}.\nonumber
\end{eqnarray}

As $H \operatorname{rank} \tilde{G}=1$ then by Theorem~\ref{th4} for $r\to
\infty$ the limit distribution of the random variable
\[
\frac{M_{r} \{ T_{n} \}-\mathbf{E}M_{r} \{ T_{n} \}
}{\sqrt{\operatorname{ \mathbf{Var}} M_{r} \{ T_{n} \}}}
\]
is the same as that of
\[
\frac{C_{(1,0,\ldots,0)} \varsigma_{1}(r)+\cdots+C_{(0,\ldots,0,1)}
\varsigma_{n+1}(r)}{\sqrt{ \operatorname{\mathbf{Var}}
(C_{(1,0,\ldots,0)} \varsigma_{1}(r)+\cdots+C_{(0,\ldots,0,1)} \varsigma
_{n+1}(r))}}=\frac{\varsigma_{1}(r)}{\sqrt{ \operatorname{\mathbf{Var}} \varsigma_{1}(r)}},
\]
where
\[
\varsigma_j(r)=\int_{\Delta(r)}H_1
\bigl(\eta_j(x)\bigr) \,\mathrm{d}x=\int_{\Delta(r)}\eta
_j(x) \,\mathrm{d}x.
\]
By Theorem~\ref{th5} the random variable ${\varsigma_{1}(r)}/{\sqrt{
\operatorname{\mathbf{Var}} \varsigma_{1}(r)}}$ is asymptotically normal with zero
mean and unit variance. By Theorem~\ref{th5} and Lemma~\ref{finint} we get
$\lim_{r\to\infty}{\operatorname{\mathbf{Var}} \varsigma_j(r)}/\break {r^{2d-\alpha}L(r)}
=c_2(d,\alpha) c_3(1,d,\alpha)$.
Finally, the application of Remark~\ref{varK} concludes the proof of
the theorem.
\end{pf*}

\begin{pf*}{Proof of Theorem~\ref{th7}}
For the function ${G} ( \cdot )$ given by (\ref{fG})
coefficients $C_{\nu}= 0$ when $\nu\in N_{1}$ or $\nu\in
N_{2}\setminus\{\nu\dvt  \mbox{exactly one } k_j=2\}$. For $\nu\in N_{2}$
with $k_j=2$ for some $j\in\{1,\ldots,m\}$, $m\ge2$, all $C_\nu$ are
equal and given below
\begin{eqnarray*}
C_{\nu}&=&\int_{\mathbb{R}^{n+m}}{G}(w)e_{(2,0,\ldots,0)}(w)\phi
\bigl(\llVert w\rrVert \bigr)\,\mathrm{d}w\\
&=&\int_{\mathbb{R}^{n+m}}\chi \biggl(
\frac{({1}/{m}) (
w_{1}^{2}+\cdots+w_{m}^{2} ) }{({1}/{n}) (w_{m+1}^{2}+\cdots
+w_{m+n}^{2} )}>a \biggr)\bigl(w_1^2-1\bigr)\prod
_{j=1}^{n+m} \frac{\mathrm{e}^{-{w_{j}^{2}}/{2}}}{\sqrt {2\pi}}\,\mathrm{d}w_j\\
&&{}+
\biggl(I_{{ma}/{(n+ma)}} \biggl(\frac{m}{2},\frac{n}{2} \biggr)-1
\biggr) \int_{\mathbb{R}}\bigl(w_1^2-1
\bigr)\frac{\mathrm{e}^{-{w_{1}^{2}}/{2}}}{\sqrt{2\pi}} \,\mathrm{d}w_1 \biggl(\int_{\mathbb{R}} \frac{\mathrm{e}^{-{w_{2}^{2}}/{2}}}{\sqrt{2\pi
}}\,\mathrm{d}w_2
\biggr)^{n+m-1} \\
&=&\frac{2 \pi^{n/2}}{(2\pi)^{(m+n)/2}\Gamma(n/2)}\frac{2 \pi
^{(m-1)/2}}{\Gamma((m-1)/2)}\\
&&{}\times\int_{\mathbb{R}}\bigl(w_{1}^2-1
\bigr) \mathrm{e}^{-{w_{1}^{2}}/{2}}\int_0^\infty
\rho^{m-2} \mathrm{e}^{-{\rho^{2}}/{2}}\int_0^{\sqrt{
{n(w_1^2+\rho^2)}/{(ma)}}}
\rho_1^{n-1} \mathrm{e}^{-{\rho_1^{2}}/{2}} \,\mathrm{d}\rho _1 \,\mathrm{d}
\rho \,\mathrm{d}w_1\\
&=& \frac{2c_4(a,n,m)}{m}.
\end{eqnarray*}

It is easy to check that for $m=1$ the above result is valid too, that
is, $C_{\nu}=2c_4(a,n,1)$.

For $\nu\in N_{2}$ with $k_j=2$ for some $j\in\{m+1,\ldots,m+n\}$ all
$C_\nu$ are equal to
\begin{eqnarray*}
C_{\nu}&=&\int_{\mathbb{R}^{n+m}}{G}(w)e_{(0,\ldots,0,2)}(w)\phi
\bigl(\llVert w\rrVert \bigr)\,\mathrm{d}w\\
&=&\int_{\mathbb{R}^{n+m}}\chi \biggl(
\frac{({1}/{m}) (
w_{1}^{2}+\cdots+w_{m}^{2} ) }{({1}/{n}) (w_{m+1}^{2}+\cdots
+w_{m+n}^{2} )}>a \biggr)
\bigl(w_{m+n}^2-1\bigr)\prod
_{j=1}^{n+m} \frac{\mathrm{e}^{-{w_{j}^{2}}/{2}}}{\sqrt{2\pi}}
\,\mathrm{d}w_j\\
&&{}+
\int_{\mathbb{R}}\bigl(w_{m+n}^2-1\bigr)
\frac{\mathrm{e}^{-{w_{m+n}^{2}}/{2}}}{\sqrt {2\pi}}
 \,\mathrm{d}w_{m+n} \biggl(\int_{\mathbb{R}}
\frac{\mathrm{e}^{-{w_{1}^{2}}/{2}}}{\sqrt{2\pi}}\,\mathrm{d}w_1 \biggr)^{n+m-1}
 \biggl(I_{{ma}/{(n+ma)}} \biggl({\frac{m}{2},\frac{n}{2}}
\biggr)-1 \biggr)\\
&=&\int_{\mathbb{R}^{n+m}} \biggl(1-\chi \biggl(
\frac{
({1}/{n}) (w_{m+1}^{2}+\cdots
+w_{m+n}^{2} )}{({1}/{m}) (
w_{1}^{2}+\cdots+w_{m}^{2} ) }>\frac{1}{a} \biggr) \biggr)
\bigl(w_{m+n}^2-1
\bigr)
\prod_{j=1}^{n+m}
\frac{\mathrm{e}^{-{w_{j}^{2}}/{2}}}{\sqrt{2\pi}}\,\mathrm{d}w_j \\
&=&-\frac{2c_4 ({1}/{a},m,n )}{n}\\
&=&-\frac{2c_4(a,n,m)}{n}.
\end{eqnarray*}

As $H \operatorname{rank} {G}=2$ then by Theorem~\ref{th4} for $r\to\infty$
the limit distribution of the random variable
\[
\frac{M_{r} \{ F_{m,n} \}-\mathbf{E}M_{r} \{ F_{m,n}
\}}{\sqrt{ \operatorname{\mathbf{Var}} M_{r} \{ F_{m,n} \}}}
\]
is the same as that of
\begin{eqnarray*}
&&\frac{C_{(2,0,\ldots,0)} \tilde{\varsigma}_{1}(r)+\cdots
+C_{(0,\ldots,0,2)} \tilde{\varsigma}_{n+m}(r)}{\sqrt{ \operatorname{\mathbf{Var}}
(C_{(2,0,\ldots,0)} \tilde{\varsigma}_{1}(r)+\cdots+C_{(0,\ldots,0,2)}
\tilde{\varsigma}_{n+m}(r))}}
\\
&&\qquad=\frac{({1}/{m})( (\tilde{\varsigma}_{1}(r)+\cdots+\tilde{\varsigma
}_{m}(r) )-({1}/{n})( (\tilde{\varsigma}_{m+1}(r)+\cdots+\tilde
{\varsigma}_{m+n}(r) )}{\sqrt{ \operatorname{\mathbf{Var}}
(({1}/{m})(
(\tilde{\varsigma}_{1}(r)+\cdots+\tilde{\varsigma}_{m}(r) )-
({1}/{n})( (\tilde{\varsigma}_{m+1}(r)+\cdots
+\tilde{\varsigma}_{m+n}(r) ) )}},
\end{eqnarray*}
where
\[
\tilde{\varsigma}_j(r)=\int_{\Delta(r)}H_2
\bigl(\eta_j(x)\bigr) \,\mathrm{d}x=\int_{\Delta
(r)}
\eta_j^2(x) \,\mathrm{d}x-\bigl|\Delta(r)\bigr|.
\]
By Theorem~\ref{th5}, we deduce that for $r\to\infty$ the
distributions of $\tilde{\varsigma}_{j}(r)/r^{d-\alpha}L(r)$ converge
to the distributions of $X_{2,j}$, where $X_{2,j}$ are independent
copies of $X_{2}$.
The application of Remark~\ref{varK} concludes the proof of the theorem.
\end{pf*}

\begin{pf*}{Proof of Theorem~\ref{th8}} It is sufficient to
investigate the case $a(r)>0$.
First, we verify condition~(\ref{con1}) for the function
%
%
\begin{equation}
\label{tG1} \tilde{G}_r ( w )=\chi \biggl( \frac{w_{1}}
{\sqrt{
({1}/{n}) (w_{2}^{2}+\cdots
+w_{n+1}^{2} )}}>a(r)
\biggr)-\frac{1}{2} I_{
{n}/{(n+a(r)^2)}} \biggl(\frac{n}{2},
\frac{1}{2} \biggr).
\end{equation}

By (\ref{par}) and (\ref{C10}) it is enough to check that
%
%
\begin{equation}
\label{tG2} \bigl(n+a(r)^2 \bigr)^{n}\int
_{R^{n+1}}\tilde{G}^2_{r}(w) \phi\bigl(
\llVert w\rrVert \bigr) \,\mathrm{d}w=\mathrm{o} \bigl(r^{\gamma/2} \bigr),\qquad r\to\infty.
\end{equation}

It follows from (\ref{tG1}) that
\begin{eqnarray*}
\int_{R^{n+1}}\tilde{G}^2_{r}(w) \phi
\bigl(\llVert w\rrVert \bigr) \,\mathrm{d}w&=&
\int_{\mathbb{R}^{n+1}}\chi \biggl(
\frac{w_{1}}{\sqrt{({1}/{n})
(w_{2}^{2}+\cdots
+w_{n+1}^{2} )}}>a(r) \biggr)\prod_{j=1}^{n+1}
\frac{\mathrm{e}^{-
{w_{j}^{2}}/{2}}}{\sqrt{2\pi}}\,\mathrm{d}w_j
\\
&&{}\times \biggl(1-I_{{n}/{(n+a(r)^2)}}
 \biggl(\frac{n}{2},\frac{1}{2}
\biggr) \biggr)+ \frac14 I^2_{{n}/{(n+a(r)^2)}}
\biggl(\frac{n}{2},
\frac
{1}{2} \biggr).
\end{eqnarray*}
For the incomplete beta function, we get
\[
I_{{n}/{(n+a(r)^2)}} \biggl(\frac{n}{2},\frac{1}{2} \biggr)=
\frac{\Gamma
((n+1)/2)}{\Gamma(n/2) \Gamma(1/2)}\int_{0}^{{n}/{(n+a(r)^2)}}
\frac
{t^{{n}/{2}-1}}{\sqrt{1-t}}\,\mathrm{d}t=\mathcal{O} \bigl( \bigl(n+a(r)^2
\bigr)^{-n/2} \bigr).
\]

Using the upper bound (7) in \cite{cook} for the complementary
cumulative distribution function of the standard normal distribution,
we conclude
\begin{eqnarray*}
&&\int_{\mathbb{R}^{n+1}}\chi \biggl( \frac{w_{1}}{\sqrt{({1}/{n})
(w_{2}^{2}+\cdots
+w_{n+1}^{2} )}}>a(r) \biggr)\prod
_{j=1}^{n+1} \frac{\mathrm{e}^{-
{w_{j}^{2}}/{2}}}{\sqrt{2\pi}}\,\mathrm{d}w_j
\\
&&\quad=\frac{2 \pi^{n/2}}{(2\pi)^{(n+1)/2}\Gamma(n/2)}s\int_0^\infty\rho^{n-1}
\mathrm{e}^{-{\rho^{2}}/{2}}\int_{a(r)\rho
/\sqrt{n}}^\infty \mathrm{e}^{-{w_{1}^{2}}/{2}}
\,\mathrm{d}w_1 \,\mathrm{d}\rho\\
&&\quad\le \frac{2 \sqrt{2}}{2^{n/2}\sqrt{\pi}\Gamma(n/2)}
\int_0^\infty\rho^{n-1}
\mathrm{e}^{-{\rho^{2}}/{2}}\frac{\mathrm{e}^{-
{a^2(r)\rho^{2}}/{(2n)}}}{{a(r)\rho}/{\sqrt{n}}+\sqrt{8/\pi+
{a^2(r)\rho^2}/{{n}}}} \,\mathrm{d}\rho \\
&&\quad= \mathcal{O} \bigl(
a^{-n}(r) \bigr),\qquad r\to\infty.
\end{eqnarray*}

Therefore,
\[
\bigl(n+a(r)^2 \bigr)^{n}\int_{R^{n+1}}
\tilde{G}^2_{r}(w) \phi\bigl(\llVert w\rrVert \bigr) \,\mathrm{d}w=
\mathcal{O} \bigl( a^{n}(r) \bigr),\qquad r\to\infty,
\]
and condition (\ref{tG2}) holds if $a(r)=\mathrm{o}(r^{\gamma/2n})$, when $r\to
\infty$.
The application of Theorems~\ref{th4} and~\ref{th5} yields the
statement of the theorem.
\end{pf*}

\begin{pf*}{Proof of Theorem~\ref{th9}}
By (\ref{fG}), we obtain
\begin{eqnarray*}
\int_{R^{n+m}}{G}^2_{r}(w) \phi\bigl(
\llVert w\rrVert \bigr) \,\mathrm{d}w&=&\int_{\mathbb
{R}^{n+m}}\chi \biggl(
\frac{({1}/{m}) (
w_{1}^{2}+\cdots+w_{m}^{2} ) }{({1}/{n}) (w_{m+1}^{2}+\cdots
+w_{m+n}^{2} )}>a(r) \biggr)\prod_{j=1}^{n+m}
\frac{\mathrm{e}^{-
{w_{j}^{2}}/{2}}}{\sqrt{2\pi}}\,\mathrm{d}w_j
\\
&&{}\times \biggl(2I_{{ma(r)}/{(n+ma(r))}} \biggl(\frac{m}{2},\frac
{n}{2}
\biggr)-1 \biggr)\\
&&{}+ \biggl(I_{{ma(r)}/{(n+ma(r))}} \biggl(\frac
{m}{2},
\frac{n}{2} \biggr)-1 \biggr)^2.
\end{eqnarray*}

The integral can be estimated as follows
\begin{eqnarray*}
&&\int_{\mathbb{R}^{n+m}}\chi \biggl( \frac{{1}/{m} (
w_{1}^{2}+\cdots+w_{m}^{2} ) }{{1}/{n} (w_{m+1}^{2}+\cdots
+w_{m+n}^{2} )}>a(r) \biggr)\prod
_{j=1}^{n+m} \frac{\mathrm{e}^{-
{w_{j}^{2}}/{2}}}{\sqrt{2\pi}}\,\mathrm{d}w_j
\\
&&\quad=\frac{4 \pi^{(n+m)/2}}{\Gamma(n/2)\Gamma(m/2)}
\frac{1}{(2\pi)^{(n+m)/2}}\int_0^\infty
\rho^{m-1} \mathrm{e}^{-{\rho
^{2}}/{2}}\int_0^{\sqrt{{n}/{(ma(r))}}\rho}
\rho_1^{n-1}\mathrm{e}^{-{\rho
_{1}^{2}}/{2}} \,\mathrm{d}\rho_1 \,\mathrm{d}\rho
\\
&&\quad\le \biggl(\frac{n}{ma(r)} \biggr)^{n/2}\frac{2^{2-(n+m)/2} }{n\Gamma(n/2)\Gamma(m/2)} \int_0^\infty\rho
^{n+m-1} \mathrm{e}^{-{\rho^{2}}/{2}} \,\mathrm{d}\rho \\
&&\quad= \mathcal{O} \bigl(
a^{-n/2}(r) \bigr),\qquad r\to\infty.
\end{eqnarray*}

By properties of the incomplete beta function, we get
\[
1-I_{{ma(r)}/{(n+ma(r))}} \biggl(\frac{m}{2},\frac{n}{2}
\biggr)=I_{
{n}/{(n+ma(r))}} \biggl(\frac{n}{2},\frac{m}{2} \biggr)=
\mathcal{O} \bigl( \bigl(n+ma(r) \bigr)^{-n/2} \bigr),\qquad r\to\infty.
\]

Therefore by (\ref{par})
\[
{\sum_{l\geq3}\sum_{\nu\in
N_{l}}
\frac{C_\nu^2(r)}{\nu!}} \Big/{\sum_{\nu\in N_{2}}\frac{C_{\nu
}^{2}(r)}{\nu!}}
\le\frac{C}{c_4^2(a(r),n,m)} \int_{R^{n+m}}{G}^2_{r}(w)
\phi\bigl(\llVert w\rrVert \bigr) \,\mathrm{d}w=\mathcal{O} \bigl( a^{n/2}(r)
\bigr),
\]
and condition (\ref{con1}) holds if $a(r)=\mathrm{o}(r^{\gamma/n})$, when $r\to
\infty$.
Steps similar to the proof of Theorem~\ref{th7} yield the statement of
the theorem.
\end{pf*}

\begin{pf*}{Proof of Theorem~\ref{th10}}
Let $G_1(\eta(x)):=\chi ( T_{n}(x)>a )$. By Assumption~\ref{ass3} we obtain
\[
M_r\{T_n\}=\int_{\Delta(r)}\chi \bigl(
T_{n}(x)>a \bigr) \,\mathrm{d}x=\int_{\Delta(r)}G_1
\bigl(\eta(x)\bigr) \,\mathrm{d}x=\int_{\Delta(r)}{\hat G_1}
\bigl(\tilde {\eta}(x) \bigr) \,\mathrm{d}x,
\]
where ${\hat G_1} (w )=G_1 (\mathcal{A}^{1/2}w )$.
By (\ref{expg}) and the orthogonality of $\mathcal{A}^{1/2}$ , we get
\begin{eqnarray*}
\mathbf{E} M_r\{T_n\}&=&\bigl|\Delta(r)\bigr|\int
_{\mathbb{R}^{n+1}}{G}_1\bigl(\mathcal {A}^{1/2}w
\bigr)\phi\bigl(\llVert w\rrVert \bigr)\,\mathrm{d}w=\bigl|\Delta(r)\bigr|\int_{\mathbb
{R}^{n+1}}{G}_1(w)
\phi\bigl(\llVert w\rrVert \bigr)\,\mathrm{d}w
\\
&=&\llvert \Delta\rrvert r^{d} \biggl(\frac12-\frac12
\biggl(1-I_{
{n}/{(n+a^2)}} \biggl(\frac{n}{2},\frac{1}{2} \biggr)
\biggr)\cdot\operatorname{sgn}(a) \biggr).
\end{eqnarray*}

$\mathcal{A}^{1/2}w$ is a linear transformation of $w$. Hence, for the
function $\tilde{G} ( \cdot )$ given by (\ref{tG}) $H \operatorname
{rank}  \tilde{G}(\mathcal{A}^{1/2}w)=$ $H \operatorname{rank} \tilde
{G}(w)=1$ and to obtain the limit theorem we need only to find the
coefficients $C_{\nu}$, $\nu\in N_{1}$, of the function $\tilde
{G}(\mathcal{A}^{1/2}w)$.\vadjust{\goodbreak}

Due to the orthogonality of $\mathcal{A}^{-1/2}$, it follows that $\sum_{i=1}^{n+1}a^2_{ji}=1$. Therefore, for $\nu\in N_{1}$ such that
$k_j=1$, by (\ref{sumH}) we obtain that
\begin{eqnarray*}
C_{\nu}&=&\int_{\mathbb{R}^{n+1}}\tilde{G}\bigl(
\mathcal{A}^{1/2}w\bigr)e_{\nu
}(w)\phi\bigl(\llVert w\rrVert
\bigr)\,\mathrm{d}w\\
&=&\int_{\mathbb{R}^{n+1}}\tilde {G}(w)e_{\nu}\bigl(
\mathcal{A}^{-1/2}w\bigr)\phi\bigl(\llVert w\rrVert \bigr)\,\mathrm{d}w
\\
&= &\int_{\mathbb{R}^{n+1}}\tilde{G}(w) \sum_{i=1}^{n+1}a_{ji}H_{1}(w_i)
\phi\bigl(\llVert w\rrVert \bigr)\,\mathrm{d}w=\frac
{a_{j1}}{\sqrt{2\pi} (1+a^2/n )^{n/2}}.
\end{eqnarray*}

Hence, for $r\to\infty$ and $\varsigma_j(r)$ defined in Theorem~\ref
{th6} the asymptotic distributions of the random variables
\[
\frac{M_{r} \{ T_{n} \}-\mathbf{E}M_{r} \{ T_{n} \}
}{\sqrt{ \operatorname{\mathbf{Var}} M_{r} \{ T_{n} \}}} \quad\mbox{and}\quad \frac{\sum_{j=1}^{n+1}a_{j1}\varsigma
_{j}(r)}{\sqrt{ \operatorname{\mathbf{Var}}
(\sum_{j=1}^{n+1}a_{j1}\varsigma_{j}(r) )}}
\]
coincide.
Note that $\sum_{j=1}^{n+1}a^2_{j1}=1$. Then, similarly to the proof of
Theorem~\ref{th6}, we get the statement of the theorem.
\end{pf*}

\begin{pf*}{Proof of Theorem~\ref{th11}} Similar to Theorem~\ref
{th10} it is easy to show that
\[
\mathbf{E} M_r\{F_{m,n}\}=\llvert \Delta\rrvert
r^{d} \biggl(1-I_{
{ma}/{(n+ma)}} \biggl(\frac{m}{2},
\frac{n}{2} \biggr) \biggr).
\]

For the function $G (\cdot )$ given by (\ref{fG}) $H \operatorname
{rank} { G}(\mathcal{A}^{1/2}w)=$ $H \operatorname{rank} { G}(w)=2$ and to
obtain the limit theorem we need only to find the coefficients $C_{\nu
}$, $\nu\in N_{2}$, of the function $G(\mathcal{A}^{1/2}w)$.

By (\ref{sumH}) and the orthogonality of both $\mathcal{A}_1$ and
$\mathcal{A}_2$, for $\nu\in N_{2}$ such that $k_j=2$, we obtain
\begin{eqnarray*}
C_{\nu}&=&\int_{\mathbb{R}^{n+m}}{G}(w)e_{\nu}\bigl(
\mathcal{A}^{-1/2}w\bigr)\phi \bigl(\llVert w\rrVert \bigr)\,\mathrm{d}w\\
&= &\int
_{\mathbb{R}^{n+m}}{G}(w) \sum_{i=1}^{n+m}a^2_{ji}H_{2}(w_i)
\phi\bigl(\llVert w\rrVert \bigr)\,\mathrm{d}w
\\
&=&2c_4(a,n,m) \Biggl( \frac{1}{m}\sum
_{i=1}^{m}a^2_{ji} -
\frac{1}{n}\sum_{i=m+1}^{m+n}a^2_{ji}
\Biggr)\\
&=& 2c_4(a,n,m) \cdot \cases{\displaystyle \frac{1}{m} ,&\quad $\mbox{if
$1\le j\le m$,}$ \vspace*{2pt}
\cr
\displaystyle-\frac{1}{n}, & \quad$\mbox{if $m+1\le j\le
m+n$,}$}
\end{eqnarray*}
while for $\nu\in N_{2}$ such that $k_j=k_l=1$, $1\le j<l\le m+n$:
\begin{eqnarray*}
C_{\nu}&=& \int_{\mathbb{R}^{n+m}}{G}(w) \sum
_{i=1}^{n+m}a_{ji}a_{li}H^2_{1}(w_i)
\phi\bigl(\llVert w\rrVert \bigr)\,\mathrm{d}w\\
&=& \sum_{i=1}^{n+m}a_{ji}a_{li}
\int_{\mathbb{R}^{n+m}}{G}(w) \bigl(H_2(w_i)+1
\bigr)\phi\bigl(\llVert w\rrVert \bigr)\,\mathrm{d}w\\
&=&2c_4(a,n,m) \Biggl(
\frac{1}{m}\sum_{i=1}^{m}a_{ji}a_{li}
-\frac{1}{n}\sum_{i=m+1}^{m+n}a_{ji}a_{li}
\Biggr)=0.
\end{eqnarray*}

%
The rest of the proof is omitted as it follows from virtually identical
arguments as in Theorem~\ref{th7}.
\end{pf*}

\section{Simulation results}\label{sec9}
To show different types of the limit behaviour for weakly and strongly
dependent models we present a simulation result based on the
theoretical findings.

For $d=2$, we chose two models of $\eta(x)$: short-range dependent
normal scale mixture field with the covariance function $\mathbf
{{B}}(\llVert  x\rrVert )=\mathcal{I}\cdot\exp (-\llVert  x\rrVert ^2 )$ and long-range dependent Cauchy field which covariance
function is $\mathbf{{B}}(\llVert  x\rrVert )=\mathcal{I}\cdot
(1+\llVert  x\rrVert ^2 )^{-1/4}$, consults \cite{sch}. We
used three independent copies of $\eta_1(x)$ to produce Fisher--Snedecor
fields $F_{1,2}(x)$, $x\in\mathbb{R}^{2}$, for each above model. The
first row of Figure~\ref{fig1} shows excursion sets above level 1 for
realizations of these two Fisher--Snedecor fields (from left to right).
The excursion sets are shown in black colour. Images in each column of
Figure~\ref{fig1} correspond to the same model. The figure was
generated by the R package \textsc{RandomFields} \cite{sch}.

Further, we simulated 1000 realizations of each $F_{1,2}(x)$ field and
computed areas of the excursion set for each realisation. Applying the
transformations given in Theorems~\ref{th3} and \ref{th7} we compared
empirical distributions of the areas to the normal law.
The second row of Figure~\ref{fig1} demonstrates normal Q--Q plots of
1000 realisations of the area of the excursion set. The observation
window was chosen to be large enough to obtain results close to the
asymptotic ones. The Q--Q plots clearly manifest differences in two
types of limit behaviour and support our findings.

\section{Conclusions}\label{sec10}
We have obtained limit distributions of the first Minkowski functional
of both weakly and strongly dependent vector random fields. In
particular, special attention was devoted to Student and
Fisher--Snedecor random fields. The techniques developed in Sections~\ref{sec4} and \ref{sec6} may be applied to other problems, which deal with
limit distributions of various functionals of vector random fields. The
analysis and the approach to the first Minkowski functional based on
functions of vector random fields are new and contribute to the
investigations of excursion sets in the former literature.

The results presented in the paper pose new problems and provide the
theoretical framework for studying more complex models.
It would be interesting:
\begin{itemize}
\item to obtain similar results for other Minkowski functionals;
\item to derive analogous results under different long-range
assumptions on covariance functions of vector random fields, consult
\cite{arc1,arc2,har};
\item to study the rate of convergence to the limit distributions,
consult \cite{leo1}.
\end{itemize}

\section*{Acknowledgements} Nikolai Leonenko was partially supported by
the grant of the Commission of the European Communities
PIRSES-GA-2008-230804 (Marie Curie) ``Multi-parameter Multi-fractional
Brownian Motion.''

The authors are grateful to the referees and Editor-in-Chief for
comments and suggestions which led to improvements in the style of the paper.


%



\printhistory

\end{document}